\documentclass[11pt]{article}
\usepackage{amsmath,amssymb} 

\begin{document}

\title{Persistently laminar branched surfaces}

\author{Ying-Qing Wu} 
\date{}
\maketitle



\begin{abstract}
  We define sink marks for branched complexes and find conditions for
  them to determine a branched surface structure.  These will be used to
  construct branched surfaces in knot and tangle complements.  We will
  extend Delman's theorem and prove that a non 2-bridge Montesinos
  knot $K$ has a persistently laminar branched surface unless it is
  equivalent to $K(1/2q_1,\, 1/q_2,\, 1/q_3,\, -1)$ for some positive
  integers $q_i$.  In most cases these branched surfaces are genuine,
  in which case $K$ admits no atoroidal Seifert fibered surgery.  It
  will also be shown that there are many persistently laminar tangles.
\end{abstract}

\newcommand{\proof}{\noindent {\bf Proof.} }
\newcommand{\qed}{\quad $\Box$}
\newtheorem{thm}{Theorem}[section]
\newtheorem{prop}[thm]{Proposition} 
\newtheorem{lemma}[thm]{Lemma} 
\newtheorem{cor}[thm]{Corollary} 
\newtheorem{defn}[thm]{Definition} 
\newtheorem{convention}[thm]{Convention} 
\newtheorem{notation}[thm]{Notation} 
\newtheorem{qtn}[thm]{Question} 
\newtheorem{example}[thm]{Example} 
\newtheorem{remark}[thm]{Remark} 
\newtheorem{conj}[thm]{Conjecture} 
\newtheorem{prob}[thm]{Problem} 
\newtheorem{rem}[thm]{Remark} 

\newcommand{\bdd}{\partial}
\newcommand{\Int}{{\rm Int}}
\renewcommand{\a}{\alpha}
\renewcommand{\b}{\beta}

\input epsf.tex

\section{Introduction}

Essential laminations are important tools in the study of topology of
3-manifolds and exceptional Dehn surgery.  Denote by $K(r)$ the
manifold obtained by $r$ surgery on a hyperbolic knot $K$ in a
closed 3-manifold.  Then the surgery is {\it exceptional\/} if $K(r)$
is non-hyperbolic, i.e.\ it is reducible, toroidal or Seifert
fibered.  A 3-manifold is {\it laminar\/} if it contains an essential
lamination.  If $M$ is laminar then it is irreducible, and if the
lamination is genuine in the sense that some complementary region is
not an $I$-bundle then $M$ is not a small Seifert fibered manifold.
In certain cases essential lamination can also be used to detect
toroidal manifolds, see for example [Wu2].

A lamination is essential if and only if it is carried by an essential
branched surface [GO].  In [Li] Li defined laminar branched surfaces.
These are essential branched surfaces that satisfy some mild extra
conditions, which makes it easier to detect such branched surfaces.
Denote by $K(r)$ the manifold obtained by Dehn surgery on a knot $K$
along slope $r$.  A laminar branched surface $\Sigma$ in the exterior
of $K$ is {\it persistently laminar\/} if it remains laminar in $K(r)$
for all non-meridional slopes $r$.

Combinatorially a branched surface is a 2-complex $\Sigma$ whose
singular set $b(\Sigma)$ is a set of immersed curves on $\Sigma$,
called the branch curves or branch loci of $\Sigma$.  A branched
surface structure on $\Sigma$ is a smooth structure in a neighborhood
of $b(\Sigma)$ so that every point has a neighborhood modeled on that
in Figure 2.1(a).  Li used sink directions on $b(\Sigma)$ to indicate
the smooth structure near $b(\Sigma)$.  See Section 2.  This is a very
useful way to define branched surface structure on $\Sigma$, however,
it is difficult to use when the branched surface is complicated.
Since there are only three ways to smooth the complex near any branch
curve, we can indicate the smooth structure on a segment $\alpha$ of
$b(\Sigma)$ using a sink mark instead, which is either an orientation
or a diamond sign on $\alpha$.  See Section 3 for more details.  This is
particularly useful for tangle complexes $\Sigma = Q \cup P \cup D$ in
the exterior of a tangle or knot $L$, where $Q$ is the tube on the
boundary of a tubular neighborhood of $L$, $P$ is a set of punctures
surfaces with $\bdd P$ a set of meridional curves on $Q$, and $D$ a
set of surfaces with boundary on $P \cup Q$ intersecting each meridian
of $Q-P$ exactly once.  The Unique Extension Lemma (Lemma 3.7) shows
that the branched surface structure of such a complex is completely
determined by the sink marks on $\bdd D$, and these sink marks define
a branched surface structure if and only if they satisfy some simple
conditions.

We will use this result to reconstruct Delman's channel surfaces [De].
By Li's result [Li] and Lemma 3.7 it is now easy to show that these
branched surfaces are laminar, and most of them are genuine.  See
Theorem 5.3.  In Section 6 we will further explore and extend Delman's
half channel surfaces construction using sink marks.  The simple
pictures of sink marks on these branched surfaces allow us to find
various semi-allowable path with desirable properties in the
Hatcher-Thurston diagram of a rational tangle.  See Proposition 6.5.
We can then strengthen the main theorem of Delman in [De], which says
that a Montesinos knot has a persistent lamination unless it is a
pretzel knot $K = K(p_1/q_1, p_2/q_2, p_3/q_3, n)$, where $p_i = 1$ or
$q_i - 1$, $q_1$ is even, $q_i$ positive, and $n=-1$ or $-2$.  The
following theorem shows that up equivalence (i.e.\ up to taking mirror
image) we must have $p_i = 1$ and $n=-1$.

\bigskip
\noindent {\bf Theorem 6.6} {\em Let $K$ be a non 2-bridge Montesinos
  knot.  Then $K$ has a persistently laminar branched surface in its
  complement unless it is equivalent to $K(1/q_1,\, 1/q_2,\, 1/q_3,\,
  -1)$, where $q_i$ are positive integers, and $q_1$ is even.}
\bigskip

As a consequence, we see that the knot $10_{142}$ is persistently
laminar.  This is one of the 5 knots with crossing number at most 10
which were not known whether surgery on them always produce laminar
manifolds [Ga, FQ 1.2].  The construction can also be applied to more
general knots.  For example, if $F$ is a minimal Seifert surface of
two component link $L$, $\alpha$ is an arc on $F$ connecting the two
components, and $K$ is obtained by replacing $N(\alpha) \cap L$ with a
$2n$-twist tangle with $|n|\geq 2$, then $K$ is persistently laminar.
See Corollary 6.9.

We are particularly interested in determining which Montesinos knots
of length 3 have a persistently laminar branched surface which is
genuine in the sense that it has a complementary region which is not
an $I$-fiber.  The following shows that most of them do have such a
branched surface.

\bigskip
\noindent {\bf Theorem 6.7} {\em 
  Let $K$ be a Montesinos knot of length 3.  Then $K$ has a
  {\em genuine\/} persistently laminar branched surface in its
  complement unless $K$ is equivalent to $K(1/q_1,\, 1/q_2,\,
  p_3/q_3,\, n)$, such that either

  (1) $n = 0$, $q_i\geq 2$, and $p_3 = 1$; or

  (2) $n =-1$, $q_i \geq 2$, and $p_3 = 1, 2$ or $q_3 - 1$.
}
\bigskip

Exceptional surgeries on arborescent knots have all been classified
except atoroidal Seifert fibered surgeries on Montesinos knots of
length 3.  By a theorem of Brittenham [Br], if $K$ has a genuine
persistently laminar branched surface then $K(r)$ is not a small
Seifert fibered manifold for any nontrivial $r$.  Using this and the
results of [Wu3] we have the following two theorem, according to
whether $K$ is pretzel or not.  Here a Montesinos knot is a {\it
  pretzel knot\/} if it can be written as $K(1/q_1,\, 1/q_2, ...,
1/q_k, n)$ for some integers $q_i$ and $n$ with $|q_i|\geq 2$, and it
is a {\it genuine pretzel knot\/} if in addition $n=0$.  The number
$k$ is called the length of $K$.

\bigskip
\noindent {\bf Theorem 7.2} {\em 
  Let $K$ be a pretzel knot of length 3.  If $K$ admits an
  atoroidal Seifert fibered surgery, then $K$ is equivalent to
  $K(\frac 1{q_1}, \frac 1{q_2}, \frac 1{q_3}, n)$ such that either
  $n=0$ and hence $K$ is a genuine pretzel knot, or $n = -1$ and
  $q_i>0$.  In either case $q_i$ satisfy $\frac 1{|q_1|-1} + \frac
  1{|q_2|-1} + \frac 1{|q_3|-1} > 1$.
}
\bigskip

Note that $q_i$ satisfies the above inequality if and only if, up to
relabeling, $(|q_1|, |q_2|, |q_3|) = (2, |q_2|, |q_3|)$, $(3,3,|q_3|)$, or
$(3,4,5)$.

\bigskip
\noindent {\bf Theorem 7.3} {\em 
  Let $K$ be a Montesinos knot of length 3.  If $K$ is not a
  pretzel knot and $K$ admits an atoroidal Seifert fibered surgery
  $K(r)$, then $K$ is equivalent to one of the following.

  (a) $K(-2/3,\, 1/3,\, 2/5)$;

  (b) $K(-1/2,\, 1/3,\, 2/(2a + 1)\, )$ and $a \in \{3,4,5,6\}$.

  (c) $K(-1/2,\, 1/q,\, 2/5)$ for some $q\geq 3$ odd;
}
\bigskip

The construction in Section 5 can also be modified to make
persistently laminar branched surfaces in tangle spaces.  A 2-string
tangle $(B,T)$ is {\it persistently laminar\/} if $B-T$ contains a
branched surface $\Sigma$, such that if $K$ is any knot that can be
written as the union of $(B,T)$ with another tangle $(B',T')$, which
is nontrivial in the sense that a curve of slope $0$ on $\bdd B$ does
not bound a disk in $B'-T'$, then $\Sigma$ is a persistently laminar
branched surface for $K$.  Brittenham [Br] showed that the tangle
$T(1/3, -1/3)$ is persistently laminar. In [Yo] Youn proved that the
tangle $T(1/3, -1/5)$ is also persistently laminar.  A Montesinos
tangle of length 2 is homeomorphic to some $T(r_1, -r_2)$ with
$0<r_i<1$ and $r_1 + r_2 \leq 1$.  The following theorem shows that
many of these are persistently laminar tangles if both $q_i$ are odd.

\bigskip
\noindent {\bf Theorem 8.5} {\em 
  If $0<r_i = p_i/q_i <\frac 23$ and $q_i$ are odd then $T =
  T(r_1, \, -r_2)$ is persistently laminar.
}
\bigskip

Thus for example if $r_i = p_i/q_i$ and $q_i\geq 3$ are old then all
tangles of type $T(1/q_1\, -1/q_2)$ or $T(r_1, -r_1)$ are persistently
laminar, and if $r_2 \in (\frac 13, \frac 23)$ then $T(r_1,\, -r_2)$
is persistently laminar for all $r_1 \in (0,1)$.  See Example 8.6 for
more details.

\medskip

The paper is organized as follows.  Section 2 gives some basic
definitions and lemmas, including laminar, pre-laminar, and
combinatorial branched surfaces, which is a branched 2-complex
$\Sigma$ with sink directions assigned, satisfying certain
combinatorial conditions.  Proposition 2.4 shows that these conditions
determine a unique branched surface structure on $\Sigma$.  Section 3
introduces sink marks and tangle complexes, and proves the Unique
Extension Lemma 3.7, which gives an easy way to detect pre-laminar
branched surfaces among tangle complexes with sink marks.  Section 4
constructs Hatcher-Thurston branched surfaces corresponding to any
edge path in the Hatcher-Thurston diagram.  It will be used later in
the constructions of other branched surfaces.  Section 5 defines
Delman channel and Delman channel surfaces $\Sigma$ corresponding to
any allowable path $\gamma$, prove Delman's theorem that such a
branched surface in a knot complement is persistently laminar, and
show that $\Sigma$ is genuine if some vertex of $\gamma$ has corner
number at least 3.  In Section 6 we extend Delman's construction of
half channel branched surfaces and prove an existence result of
semi-allowable path with certain properties (Proposition 6.5), which
is then used to prove the existence theorems (Theorems 6.6 and 6.7)
for persistently laminar branched surfaces.  These are applied in
Section 7 to study atoroidal Seifert fibered Dehn surgery.  Section 8
constructs persitently laminar branched surfaces in tangle spaces.

\bigskip

All manifolds are assumed compact, connected and orientable unless
otherwise stated.  For any submanifold $Y$ in $X$, denote by $X|Y$ the
manifold obtained by cutting $X$ along $Y$.  When a rational number
$r$ is written as $p/q$ is is always assumed that $p,q$ are coprime.


\section{Combinatorial description of branched surfaces}

We refer the readers to [GO] for basic definitions such as essential
lamination, essential branched surface $F$, its regular neighborhood
$N(F)$, its $I$-fibers, the collapsing map $\pi: N(F) \to F$, the
horizontal surface $\bdd_h F$, and the vertical surface $\bdd_v F$.
The vertical surface $\bdd_vF$ is also called the {\it cusps\/} of
$F$.  Recall that the {\it branch locus\/} $b(F)$ of a branched
surface $F$ is the set of points which does not have a disk
neighborhood.  It is a finite union of immersed curves on $F$.  The
set of double points of $b(F)$ is denoted by $s(F)$, called the {\em
  singular points\/} of $F$, which cuts $b(F)$ into arcs and circles,
called the {\it branch curves}.  A point in $b(F) - s(F)$ has a
neighborhood which is the union of three disks $F_1, F_2, F_3$ joined
at an arc in $b(F) - s(F)$, and a point in $s(F)$ has a neighborhood
modeled on Figure 2.1(a).  Li [Li] uses an arrow to indicate the {\em
  sink direction\/} of $\alpha$.  It is an arrow pointing from
$\alpha$ into one of the surfaces attached to it, so that if the
vertical surface of a regular neighborhood of $F_1 \cup F_2 \cup F_3$
lies between $F_2$ and $F_3$ then the sink direction will point from
$\alpha$ into $F_1$, as shown in Figure 2.1(a).  We say that $\alpha$
is a {\em sink edge\/} of $F_1$ and a {\em source edge\/} of $F_2$ and
$F_3$.  We will also consider $\alpha$ as a cusp on the side of $F_2
\cup F_3$ opposite to the surface $F_1$ since it corresponds to a
piece of the cusps $\bdd_v F$ on that side.  The introduction of sink
direction greatly simplifies the way to draw branched surfaces since
now we only need to draw the topological 2-simplex and the sink
directions and do not have to draw it with the actual tangency or
smooth structure.  For example, we can use the 2-complex in Figure
2.1(b) to denote the branched surface in Figure 2.1(a).

\bigskip
\leavevmode

\centerline{\epsfbox{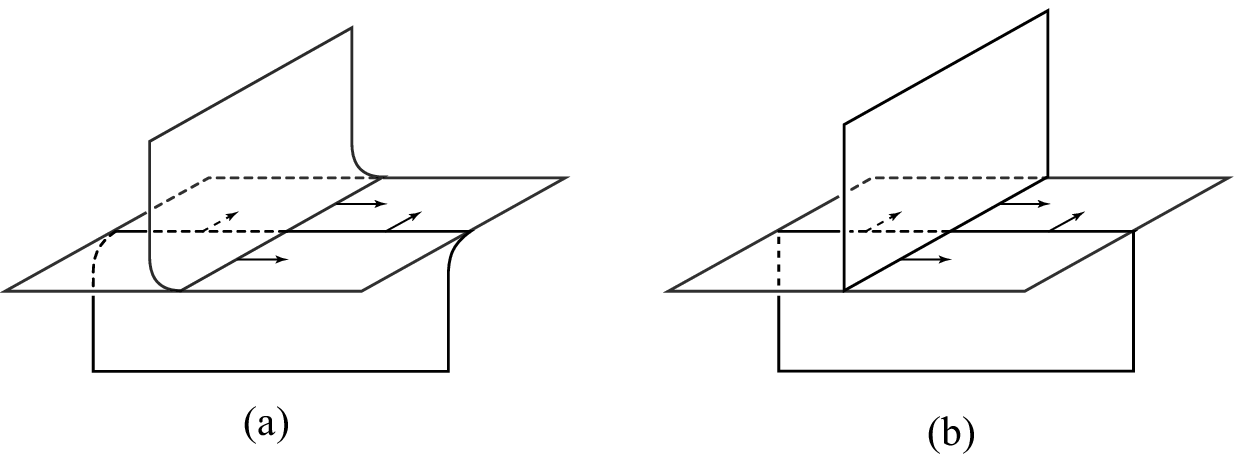}}
\bigskip
\centerline{Figure 2.1}
\bigskip

The branch locus $b(F)$ cuts $F$ into several surfaces, called the
{\em branches\/} of $F$.  We allow a branched surface $F$ to have
nonempty boundary $\bdd F$, which is a train track.  Thus the boundary
of a branch of $F$ is a union of sink arcs, source arcs, and possibly
some arcs on $\bdd F$.  A disk branch $D$ of $F$ is called a {\it sink
  disk\/} if $\bdd D$ contains some sink edges but no source edge.
This matches the definition in [Li] when $F$ has no boundary.
Similarly, a disk is a {\it source disk\/} if it has source edges but
no sink edge, and a {\it passing disk\/} if it contains
both sink edges and source edges.  $F$ is {\it sinkless\/} if it
contains no sink disk.

A {\em cusped manifold\/} is a pair $(M, \gamma)$, where $M$ is a
compact orientable 3-manifold, and $\gamma$ is a set of simple closed
curves on $\bdd M$.  Denote by $\bdd_v (M)$ a regular neighborhood of
$\gamma$ on $\bdd M$, called the vertical boundary, and by $\bdd_h(M)$
the surface $\bdd M - \Int \bdd_v(M)$, called the horizontal surface.
If $\Sigma$ is a branched surface in a 3-manifold $Y$, we use
$E(\Sigma)$ to denote $Y - \Int N(\Sigma)$, called the {\it
  exterior\/} of $\Sigma$.  It is a cusped manifold with $\gamma$ the
central curve of $\bdd_v(\Sigma)$, so $\bdd_h(M) = \bdd_h(\Sigma)$ and
$\bdd_v(M) = \bdd_v(\Sigma)$.  A disk $D$ in $M$ is a {\it monogon\/}
if $\bdd D$ intersects $\gamma$ transversely at a single point.  A
cusped manifold $(M, \gamma)$ is {\it essential\/} if $M$ is
irreducible, has no monogon, and $\bdd_h M$ is incompressible and has
no sphere component.  If $M$ is a solid torus, the {\it cusp winding
  number\/} of $(M, \gamma)$ is the minimal intersection number
between $\gamma$ and a meridian of $M$.  It is easy to see that in the
definition of essentiality of $M$, the condition that $M$ has no
monogon can be replaces by the weaker condition that no component of
$M$ is a solid torus with cusp winding number $1$.

A surface {\it carried by $\Sigma$\/} is an embedded surface in
$N(\Sigma)$ transverse to the $I$-fibers.  Let $\pi: N(\Sigma) \to
\Sigma$ be the collapsing map, which maps each $I$-fiber to a single
point of $\Sigma$.


\begin{defn} (1) A branched surface $\Sigma$ is {\em atoroidal\/} if
  any torus carried by $\Sigma$ is parallel to a component of $\bdd_h
  \Sigma$ in $N(\Sigma)$.

  (2) An {\em embedded\/} sphere $S$ in $\Sigma$ is a {\em trivial
    bubble\/} if one side of $S$ has no branch attached.  In this case
  $S$ is parallel to a spherical component $S'$ of $\bdd N(F)$ such
  that $\pi: S' \to \Sigma$ is injective.

  (3) A branched surface $\Sigma$ is {\em pre-laminar\/} if it is
  sinkless, atoroidal, and has no trivial bubble.

  (4) A closed branched surface $\Sigma$ embedded in a closed
  3-manifold $M$ is {\em laminar\/} if it is sinkless, it has no
  trivial bubble, it carries no Reeb torus of $M$, and $E(\Sigma)$ is
  an essential cusped manifold.

  (5) A branched surface $\Sigma$ in a 3-manifold $M$ with $\bdd
  \Sigma \subset \bdd M$ is {\em genuine\/} if at least one
  component $Y$ of $E(\Sigma)$ in the interior of $M$ is not an
  $I$-bundle with $\bdd_h Y$ the corresponding $\bdd I$-bundle.
\end{defn}

Li [Li, Theorem 1] proved that a laminar branched surface in a closed
orientable 3-manifold carries an essential lamination and hence is an
essential branched surface.  A laminar branched surface does not have
to be pre-laminar because it may be toroidal.  Neither is a
pre-laminar branched surface in a closed 3-manifold $M$ necessarily
laminar.  Being pre-laminar is an intrinsic property and is independent
of the embedding of $\Sigma$ in a 3-manifold, hence it can be
determined before it is embedded in $M$.  The following result follows
from [Li, Theorem 1] immediately since a pre-laminar branched surface
carries no torus and hence no Reeb torus.

\begin{lemma} If $\Sigma$ is a closed pre-laminar branched surface,
  then its embedding in a compact orientable 3-manifold $M$ is laminar
  (and hence essential) if and only if $E(\Sigma)$ is an essential
  cusped manifold. \qed
\end{lemma}

We note that the above lemma works for closed branched surface only.
To extend it to branched surfaces with boundary we need to modify the
definition of pre-laminar branched surfaces and essential cusped
manifold, for example there should be no trivial half-bubble on
$\Sigma$ and the horizontal surface of $E(\Sigma)$ must be boundary
incompressible.  In the constructions below we will need to construct
pre-laminar branched surfaces with boundary, but these will eventually
be combined to make closed branched surfaces.  Therefore we only
need the above lemma when $\Sigma$ is a closed branched surface.

Using sink directions we can define branched surface combinatorially.

\begin{defn} 
  (1) A compact 2-complex $F$ is a {\em branched complex\/} if it is
  locally modeled on the 2-complex in Figure 2.1(b) up to
  homeomorphism (without sink arrows specified).  Denote by $b(F)$ the
  set of points which does not have a surface neighborhood, and by
  $s(F)$ set of points in $b(F)$ which does not have an arc
  neighborhood in $b(F)$, called the singular set of $F$.

  (2) A {\em sink direction\/} on a component $\alpha$ of $b(F)-s(F)$
  is a vector pointing from $\alpha$ into one of the three incidented
  surfaces.

  (3) A {\em combinatorial branched surface\/} is a branched complex
  $F$ with a sink direction assigned for each component of $b(F)-s(F)$,
  such that for each point $p \in s(F)$, the 4 arcs $e_1, ..., e_4$ in
  $b(F)$ and the 6 disks $A_{ij}$ ($1\leq i<j \leq 4$) in a
  neighborhood of $p$ can be labeled so that $A_{ij} \cap b(F) = e_i
  \cup e_j$, $A_{12}$ is a sink disk, $A_{2,3}, A_{1,4}$ are passing
  disks, and the others are source disks.
\end{defn}

The following proposition shows that a combinatorial branched surface
is a branched surface after smoothing along the branch loci according
to the sink directions.  We will therefore consider any combinatorial
branched surface as a branched surface, and vise versa.

\begin{prop} Any combinatorial branched surface $F$ is homeomorphic to 
a branched surface $F'$ with sink directions preserved.  If $F$ is
embedded in a 3-manifold $M$ then it is isotopic to a branched surface
with sink directions preserved.
\end{prop}

\proof We assume $F\subset M$.  The other case is similar.  There is
no problem smoothing $F$ along $b(F) -s(F)$ according to the sink
direction, so we need only verify that the neighborhood of a point
$p\in s(F)$ can be deformed to a branched surface according to the
sink directions.  Let $D$ be the disk $A_{12} \cup A_{23} \cup A_{34}
\cup A_{14}$.  By definition $A_{12}$ is a sink disk, $A_{34}$ is a
source disk, and the other two are passing disks; hence each $e_i$
is a sink edge of one disk and a source edge of another in $D$, so up
to isotopy we may assume $D$ is a smooth disk in $M$, and $c' = e_1
\cup e_3$ and $c'' = e_2\cup e_4$ are smooth arcs on $D$.  

Note that each branch curve is a sink edge of exactly one of the three
incidented disks.  The disks $A_{13}$ intersect $D$ at $c' = e_1 \cup
e_3$.  Since the sink directions of $e_1, e_3$ points into $A_{12}$
and $A_{23}$ respectively, which are on the same side of $c'$, the
tangency of $A_{13}$ along $e_1, e_3$ matches at $p$, so $D \cup
A_{13}$ is a branched surface.  Since $A_{24}$ and $A_{13}$ has
disjoint interior while their boundaries intersect transversely at
$p$, they must lie on different sides of $D$.  For the same reason,
the sink directions on the two boundary edges $e_2, e_4$ of $c''$
points to the same side of $c''$, so the smoothing along these two
edges matches at $p$, hence $D \cup A_{13} \cup A_{24}$ is a branched
surface after smoothing.  \qed

\begin{remark} {\rm (1)  The sink directions defines a branched surface
    structure in a neighborhood $X$ of $p \in s(F)$ if and only if (i)
    there is exactly one sink disk, two passing disks, and three
    source disks, and (ii) the two passing disks intersect only at
    $p$, and the sink directions of their source edges point into the
    same disk.  Thus among the $4^3 = 81$ possible choices of sink
    directions, only $12$ of them make $X$ a branched surface. 

    (2) The sink directions of a branched surface near a singular
    point $p$ is completely determined by those of the two passing
    disks in a neighborhood of $p$.  }
\end{remark}

\begin{example} {\rm The 2-simplices in Figure 2.2(a), (b) and (c) are
    not branched surfaces.  The surface in (a) has two sink disks, the
    one in (b) has no sink disk, and in (c) the two passing disks
    have one edge in common, which implies that it is not a branched
    surface by Remark 2.5.  One can check that the 2-complex in Figure
    2.2(d) satisfies Definition 2.3(3) and hence is a branched
    surface.  It is homeomorphic to that in Figure 2.2(d), which is
    the same as the branched surface in Figure 2.1(b). }
\end{example}

\bigskip
\leavevmode

\centerline{\epsfbox{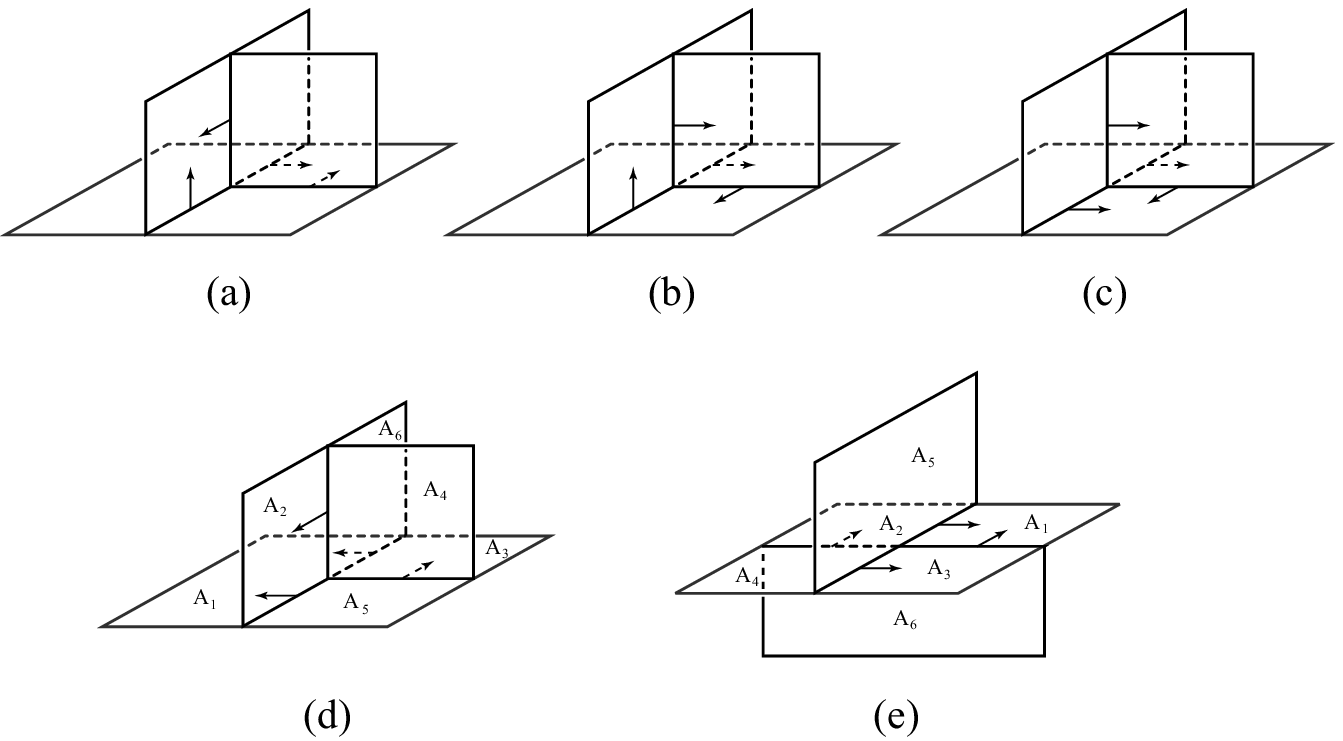}}
\bigskip
\centerline{Figure 2.2}
\bigskip

\section{Sink marks for branched surfaces}

While sink direction makes it possible to define branched surfaces
combinatorially, practically it is still very difficult using it to
define branched surface structure on branched complex with more than
just a very few branch curves.  We need to further simplify it in
order to use it to illustrate the branched surfaces to be constructed.
Since each branch curve $\a$ has three possible sink directions, we
can use an orientation of $\a$ and a diamond sign to indicate such a
choice.

More explicitly, define a {\em surface decomposition\/} of a 2-complex
$\Sigma$ to be a set of compact surfaces $S_1, ..., S_n$ in $\Sigma$
with mutually disjoint interiors, such that $\cup S_i = \Sigma$, and
each component $\alpha$ of $b(\Sigma) - s(\Sigma)$ is on the boundary
of one $S_i$ and the interior of another $S_j$.  Thus in a
neighborhood of $\alpha$, $\Sigma$ is obtained by attaching $S_i$ to
$S_j$ along the branch curve $\alpha$.

\begin{defn} 
  Suppose $\{S_i\}$ is a surface decomposition of $\Sigma$.  Then a
  {\em sink mark\/} on $\alpha \subset \bdd S_i \cap \Int S_j$ is a
  diamond sign or an orientation of $\alpha$.  The edge $\alpha$ is a
  called {\em diamond edge\/} or an {\em oriented edge\/} accordingly.
  It determines a sink direction on $\alpha$ as follows.

  (1) The sink direction of a diamond edge on $\bdd S_i$ points into
  $S_i$.

  (2) The orientation of an oriented edge $\alpha \subset \bdd S_i$
  defines a local orientation of $S_i$ and the sink direction points
  to the side of positive normal direction of $S_i$.  When $M = S^3$
  we use the orientation convention that the sink direction points to
  the right of $\alpha$ when standing on $S_j$ on the side of $S_i$
  and facing to the direction of the orientation of $\alpha$.

  (3) A set of sink marks, one for each branch curve, defines a {\em
    branched surface structure on $X$\/} if $X$, with sink directions
  determined by the sink marks, is a branched surface.
\end{defn}

For example, the 2-complexes in Figure 2.2(a)-(c) have natural
surface decompositions $(S_1, S_2, S_3)$ with each $S_i$ a flat disk.
The corresponding sink marks are given in Figure 3.1(a)-(c),
respectively.  Note that when a branched surface $\Sigma \subset S^3$
is mapped to $\Sigma'$ by an orientation reversing map of $S^3$ then
all the orientation marks of the branch curves are reversed because
the global orientation has reversed.

\bigskip
\leavevmode

\centerline{\epsfbox{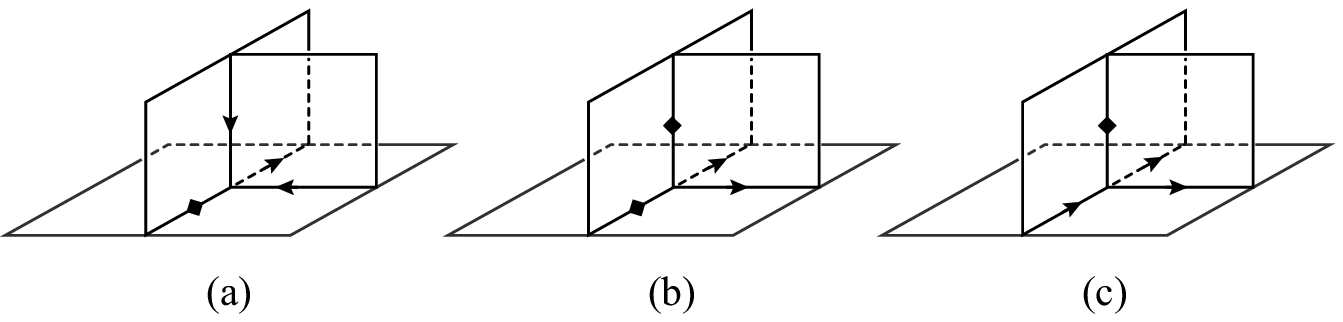}}
\bigskip
\centerline{Figure 3.1}
\bigskip

A curve $C$ with several segment marked by sink marks is {\it
consistently oriented\/} if it has an orientation which matches all
sink marks on it.  In particular, it has no diamond marks.  For
example, the central circle in Figure 3.2(3) is consistently oriented
while those in Figure 3.2(1)-(2) are not.  Also, three of the six
rectangles in Figure 3.2(1)-(3) are consistently oriented.

Now consider the branched complex $X$ shown in Figure 3.2(1).  It has an
obvious surface decomposition $(Q,P,D)$, where $Q$ is the vertical
annulus, $P$ the horizontal annulus with inner boundary attached to
$Q$, and $D = D_1 \cup D_2$ is a pair of disks, such that each of
$\alpha_i = D_i \cap P$ and $\beta_i = D_i \cap Q$ is a single arc.

\begin{defn} When both segments of $\bdd P$ are diamond edges as shown
  in Figure 3.2(2), $\bdd P$ is called a {\em meridional cusp\/} of
  $X$.  Note that in this case the cusp corresponding to $\bdd P$ is
  on the inside side of $Q$.
\end{defn}

\bigskip
\leavevmode

\centerline{\epsfbox{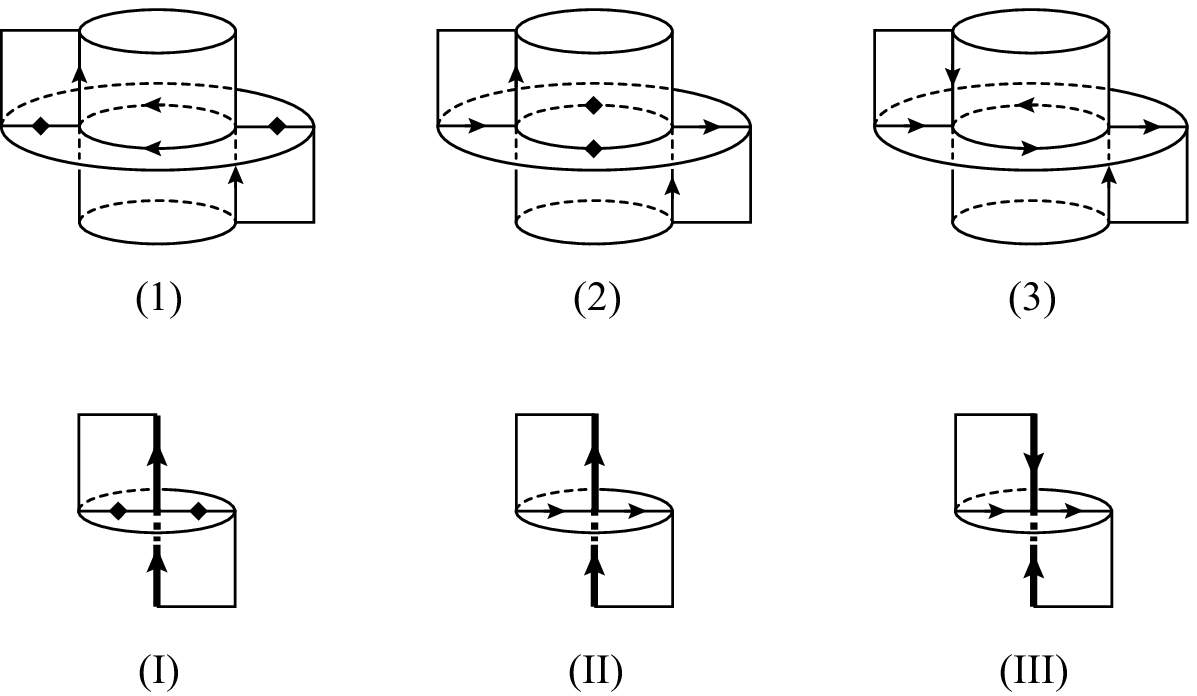}}
\bigskip
\centerline{Figure 3.2}
\bigskip

\begin{defn} Let $X$ be the underlying 2-complex in Figure 3.2(1) with
  surface decomposition $(Q,P,D)$ as above.  A set of sink marks
  assigned on the four edges of $\bdd D$ (but not on $\bdd P$) is of
  type (I) if $\bdd D \cap P$ are diamond edges and the two edges of
  $\bdd D \cap Q$ on $Q$ point to the same direction, of type (II) if
  $\bdd D$ has no diamond edge and each $\bdd D_i$ is consistently
  oriented, and of type (III) if $\bdd D$ has no diamond edge and
  exactly one $\bdd D_i$ is consistently oriented.
\end{defn}

Thus up to homeomorphism of $X$, $\bdd D$ is as shown in Figure
3.2(1)--(3) respectively, with possibly the orientations of both
segments of $\bdd D_i$ reversed for one or both $D_i$ if it is of type
(II) or (III).  Since there is no specification for sink marks on
$\bdd P$, we can shrink $Q$ and use a thick arc $K$ to represent $Q$
without loss of sink mark information, so the $X$ with sink marks in
Figure 3.2(1)-(3) can be represented by those in Figure 3.2(I), (II)
and (III), respectively.  The sink marks on $\bdd D$ induces a
piecewise orientation on $K$, called the {\it induced orientation}, or
the orientation induced by the sink marks.  The following result shows
that these sink mark systems can be uniquely extended over $\bdd P$ to
define a branched surface structure on $X$.

\begin{lemma} Let $X = Q \cup P \cup D$ be as above.  Then a sink mark
  system on $\bdd D$ with no diamond on $\bdd D \cap Q$ can be
  extended to a branched surface structure of $X$ if and only if it is
  of type (I), (II) or (III), in which case the extension is unique.
  In particular, $\bdd P$ is a meridional cusp if and only if $X$ is
  of type (II).
\end{lemma}

\proof Let $D_1 $ be the upper left disk and $D_2$ the lower right
disk of $D$ in Figure 3.2(1)-(3).  One can check that the sink
directions corresponding to the sink marks in Figure 3.2(1)-(3)
satisfy the conditions in Definition 2.3(3), hence determine a
branched surface structure in each case.  Therefore the extensions
exist.  To prove the uniqueness, assume $X$ has been assigned sink
mark system so that $X$ is a branched surface.  There are three
possibilities for the sink marks on $\bdd P$.  We want to show that
each of them corresponds to one of the types above, and the sink marks
on $\bdd P$ are completely determined by those on $\bdd D$.

{\sc Case 1.} {\it One edge of $\bdd P$ has a diamond mark.\/} In this
case this edge has a cusp on the inside of $Q$, which must extend to
the other edge on $\bdd P$ because by assumption $\bdd D\cap Q$ are
not diamond edges.  Thus both edges of $\bdd P$ are diamond edges.  In
this case $\bdd D_i$ does not pass cusp, hence must be consistently
oriented, as shown in Figure 3.2(2), possibly with orientations of
both segments of $\bdd D_i$ reversed for one or both $i$.  Therefore
$X$ is of type (II).

{\sc Case 2.} {\it The two segments of $\bdd P$ are inconsistently
  oriented.\/} In this case, near each singular point of $X$ the
branch on $Q$ with interior disjoint from $D$ is a passing disk, hence
by Remark 2.5(2), if these sink marks defines a branched surface
structure then the sink marks on $\bdd D$ are completely determined by
those on $\bdd P$.  Since the ones in Figure 3.2(1) do define a
branched surface structure, it follows that the sink marks on $\bdd D$
must be as in Figure 3.2(1) if $\bdd P$ is oriented as in the figure.
Similarly if the orientations of both segments of $\bdd P$ are
reversed then the sink marks on $\bdd D$ are obtained by reversing the
orientations of $\bdd D \cap Q$.  In either case $X$ is of type I.
Note that the orientations on $\bdd P$ are also determined by the sink
marks on $\bdd D$ as they have to point to the disk $D_i$ such that
$D_i \cap Q$ has a tail at $\bdd D_i \cap \bdd P$.

{\sc Case 3.}  {\it $\bdd P$ is consistently oriented.\/} Assume that
the orientation of $\bdd P$ is as shown in Figure 3.2(3), then the
cusp on $\bdd P$ is on the upper side of $P$.  There are two ways to
attach each $D_i$.  $\bdd D_1$ passes through the cusp and must have
inconsistent orientations on the two segments, and the orientations on
$\bdd D_2$ are consistent since it does not pass a cusp.  Therefore
$X$ is of type (III).  Similarly if the orientation of $\bdd P$ is
reversed then $\bdd D_1$ is consistently oriented while $\bdd D_2$ is
inconsistently oriented.  Hence the orientation of $\bdd P$ is
determined by the sink marks on $\bdd D$ according to the fact that
the cusp at $\bdd P$ is on the side of the disk $D_i$ whose boundary
is inconsistently oriented.

We have shown that the three possibilities above correspond to the
three types of $X$ and in each case the sink marks on $\bdd P$ are
also completely determined by those on $\bdd D$.  This completes the
proof of the lemma.  \qed
\medskip

\begin{defn} A {\em tangle complex\/} is a branched complex $X$ in
  $S^3$ with surface decomposition $(Q, P, D)$ as follows.  Suppose
  $K\subset S^3$ is a compact 1-manifold and $S$ is a set of spheres
  such that $S_i \cap K \neq \emptyset$ for any component $S_i$ of
  $S$.  Let $N(K) = K \times D^2$ be a tubular neighborhood of $K$,
  $Q$ the tubes $K \times \bdd D^2$, $P$ the punctured spheres $S -
  \Int N(K)$, and $D$ a set of compact surfaces attached to $Q\cup P$
  along some boundary curves $\gamma$ of $D$.  They satisfy (i) each
  component of $\bdd P$ is either a component of $\bdd Q$, or it has a
  regular neighborhood modeled on the underlying 2-complex in Figure
  3.2(1); (ii) each meridian curve of $Q-P$ intersects $D$ at exactly
  one point, and (iii) each component of $D$ has some boundary arcs on
  $P$ and some on $Q$.
\end{defn}

Lemma 3.4 allows us to use a thick curve to represent a tube $Q$ when
a tangle complex is a branched surface because there is no need to
specify the sink marks on $\bdd P$, as long as each point of $K \cap
S$ has a neighborhood of type (I), (II) or (III) as shown in Figure
3.2.  Thus when drawing $X$, we will simply draw $(K, S, D)$, with the
understanding that the thickened curve $K$ represents the tube $Q$
around $K$, and a disk transverse to $K$ in the figures represents a
punctured spheres $P$.

A point $p$ in the singular set $s(X)$ of a tangle complex $X = Q \cup
P \cup D$ is either on $\bdd P$ or in the interior of $P$.  In the
latter case the neighborhood of $p$ consists of one subdisk $P_0$ on
$P$ and two subdisks $D_1, D_2$ of $D$, one attached on each side of
$P_0$.  See Figure 3.3.

\bigskip
\leavevmode

\centerline{\epsfbox{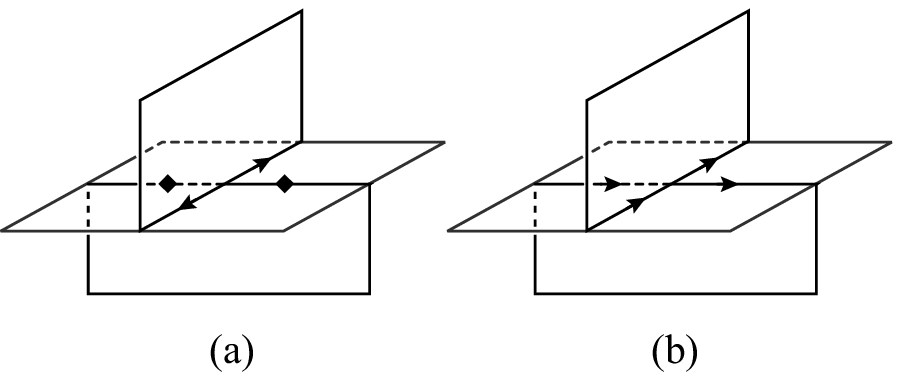}}
\bigskip
\centerline{Figure 3.3}
\bigskip

\begin{defn} Let $\Sigma = Q \cup P \cup D$ be a tangle complex with
  sink marks specified on $\bdd D$.  Then a singular point $p\in
  s(\Sigma)$ is of type (I), (II) or (III) if it lies on a component
  of $\bdd P$ which is of type (I), (II) or (III), respectively.  It
  is of type (IV) or (V), if it lies in the interior of $P$ and has a
  neighborhood as shown in Figure 3.3(a) and (b) respectively,
  possibly with orientations of both segments of $\bdd D_i \cap P$
  reversed for one or both $i$, where the horizontal disk is a subdisk
  of $P$ and the other two disks $D_i$ are in $D$.
\end{defn}

\begin{lemma} {\rm (The Unique Extension Lemma)\/}\; Let $\Sigma = Q
  \cup P \cup D$ be a tangle complex with sink marks specified on
  $\bdd D$, such that each $p \in s(\Sigma)$ is of one of the types
  (I)-(V).  Then these sink marks can be uniquely extended over $\bdd
  P$ to a branched surface structure on $\Sigma$.  Moreover, if each
  branch on $P$ has a diamond edge on its boundary then $\Sigma$ is
  pre-laminar.
\end{lemma}

\proof Checking the sink marks in Figure 3.3, one can see that
Condition (1) implies that the regular neighborhood of any singular
point $p \in s(\Sigma)$ in the interior of $P$ is a branched surface.
All other points of $s(X)$ are in $\bdd P$, and condition (2) and
Lemma 3.4 shows that there is a unique choice of sink marks for the
branch curves in $\bdd P$ to make a neighborhood of $\bdd P$ a
branched surface.  

By assumption $P$ contains no sink disk because each branch $B$ has a
diamond edge, which is a source edge for $B$.  Since each component of
$\bdd P$ is of type (I), (II) or (III), no component of $D\cap Q$ is a
diamond edge.  Since each component $D_i$ of $D$ has some boundary
edge on $Q$, which is a source edge for $D_i$, it follows that there
is no sink disk on $D$.  By definition each component $Q_i$ of $Q|P$
intersects $D$ at exactly one edge $\alpha$, hence the branch on $Q_i$
incidents $\alpha$ twice, so it has a source edge.  It follows that
$\hat \Sigma$ is sinkless.

A trivial bubble $R$ is the union of some branches.  Clearly $P \cup
Q$ contains no sphere, so $R$ must contain a component $D_i$ of $D$.
By definition $D_i$ has a boundary arc $\alpha$ on some component
$Q_j$ of $Q$ cut along $\bdd P$.  Since $Q_j$ cut along $\alpha$ is a
branch, $R$ contains $D_i \cup Q_j$, contradicting the assumption that
$R$ is a sphere.  Hence $\hat \Sigma$ contains no trivial bubble.

Suppose $S$ is a closed surface carried by $\Sigma$.  If $S$
intersects a fiber of a branch on $P$ then it must flow out of the
diamond edge of this branch into some branch $F_i$ on $D$, and since
$F_i$ has a source edge on $Q$ it must flow into $Q$.  Let $m$ be a
meridian loop of $Q$.  Then its preimage under the collapsing map
$\varphi: N(\hat \Sigma) \to \hat \Sigma$ is an annulus $A$ which is
$I$-fibered, and $\bdd A = \alpha \cup \beta$, where $\alpha$ is an
arc on a single $I$-fiber, and $\beta$ is the union of a circle and an
arc transverse to the $I$-fibers.  If $S$ is not in $N(Q)$ then by the
above it must flow into $Q$, so it intersects the $I$-fiber $\alpha$
for some meridian loop $m$.  Since $S$ is transverse to the
$I$-fibers, each component of $S\cap A$ is a curve in $A$ with a
single boundary point on $\alpha$, hence $S\cap A$ is not a compact
curve, contradicting the assumption that $S$ is a compact surface.
Therefore any connected closed surface $S$ carried by $\hat \Sigma$
must be disjoint from fibers of $P$ and $D$, so it is carried by $Q$
and hence is a torus carried by a component $Q'$ of $Q$.  Note that no
component of $\bdd P$ on $Q'$ can be marked by diamond as otherwise it
would be a source edge of a branch of $Q'$ and hence $S$ would flow
out of $Q'$ into $P$, contradicting the fact that it is carried by
$Q'$.  Therefore the inside side of $\bdd N(Q')$ is a torus $T$ on
$\bdd N(\Sigma)$ with no cusp and hence is a component of $\bdd_h
(\Sigma)$, and $S$ is parallel to $T$.  It follows that $\hat \Sigma$
is atoroidal.  \qed \medskip

\section{The Hatcher-Thurston branched surfaces}

Consider $S^3$ as $S^2 \times [-\infty, \infty]$ with each $S^2 \times
\{\pm \infty\}$ pinched to a point.  Denote by $S[x]$ the image of
$S^2 \times \{x\}$ and by $S[a,b]$ the image of $S^2 \times [a,b]$.

Let $X$ be the tangle complex in $S[a,b]$ ($a,b$ finite) with surface
decomposition $(Q, P, D)$ shown in Figure 4.1(a), where $Q$ consists
of four vertical tubes represented in the figure by 4 vertical arcs
$K$, $P$ is a horizontal punctured sphere, and $D$ the union of four
rectangles, each having two edges on $Q$, one edge on $P$ and one edge
on $\bdd S[a,b]$.  The 4 edges of $D\cap P$ are diamond edges.  Fix an
orientation of $K$ (so the two segments of any component of $K$ are
oriented consistently), which induces sink marks on $\bdd D \cap Q$.
Then we see that the neighborhood of any component of $\bdd P$ is a
complex of type (I) as in Figure 3.2(I).  Therefore by Lemma 3.7 these
sink marks extend to a unique branched surface structure on $\Sigma =
Q \cup P \cup D$, which is pre-laminar.

\bigskip
\leavevmode

\centerline{\epsfbox{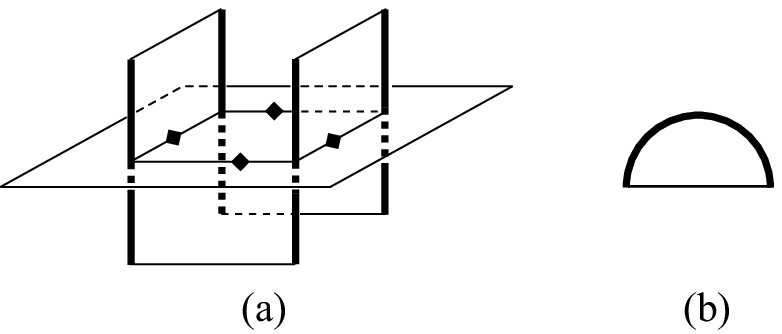}}
\bigskip
\centerline{Figure 4.1}
\bigskip

Given two rational numbers $r_i = p_i/q_i$, denote by
$\Delta(p_1/q_1,\, p_2/q_2) = |p_1 q_2 - p_2 q_1|$.  When $\Delta(r_1,
r_2) = 1$, we can deform the branched surface above by twisting the
four components of $K$ around each other so that the two top edges on
$S[0]$ have slope $r_1$ and the bottom edges have slope $r_2$.  (See
[HT] for definition of slopes of such curves.)  We will call this
branched surface the {\it Hatcher-Thurston branched surface\/}
associated to the edge $e$ from $r_1$ to $r_2$, denoted by $\Sigma(e)$
or $\Sigma(r_1, r_2)$.  This branched surface has the following
properties.  (i) When view from above, we see a pair of top edges with
slop $r_1$ on $S[a]$, and a pair of cusps on $P$ with slope $r_2$;
(ii) when view from below, we see a pair of bottom edges with slope
$r_2$ and a pair of cusps with slope $r_1$.

If $K$ is an oriented curve in a compact 3-manifold $M$, then a point
of $K$ on $\bdd M$ is {\it positive\/} if the orientation of $K$
points outward at that point, and {\it negative\/} otherwise.  The
boundary train tracks of $\Sigma(e)$ and the branched surfaces to be
constructed below depend on the local orientations of the thick arcs
$K$.  These determine the sink directions and hence the way $D$ is
attached to $Q$.  See Figure 4.2(a)-(c) for the three possible
boundary train tracks when viewed from outside of $S[a,b]$, which will
be said to be {\em positively oriented, negative oriented, and
  antiparallel}, respectively.

\bigskip
\leavevmode

\centerline{\epsfbox{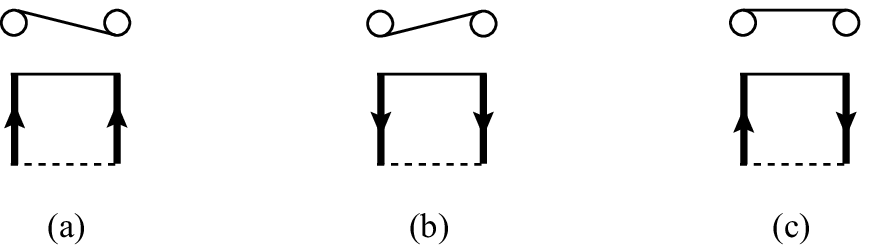}}
\bigskip
\centerline{Figure 4.2}
\bigskip

We refer the readers to [HT, Fig.\ 4] for the Hatcher-Thurston diagram
$\cal D$, which is a graph on a disk $D^2$ having $\Bbb Q \cup
\{\infty\} \subset \bdd D^2$ as vertices, with an edge connecting
$r_1$ to $r_2$ if $\Delta(r_1, r_2) = 1$.  A path on $\cal D$ is {\it
  minimal\/} if two successive edges do not lie on a triangle of $\cal
D$.  Let ${\cal D}(p/q)$ be the subdiagram of $\cal D$ consisting of
the edges of all minimal paths from $1/0$ to $p/q$.  See Figure 4.3
for ${\cal D}(3/11)$.  When $q\neq 1$ ${\cal D}(p/q)$ is a graph on a
disk $D$ containing $\bdd D$, with all vertices on $\bdd D$ and all
faces triangles.  The edges on $\bdd D$ form two paths from $1/0$ to
$p/q$, called the {\it upper boundary path\/} and the {\it lower
  boundary path\/}, with the upper one containing the vertices with
label $r_i > p/q$.  Edges on $\bdd D$ are {\it boundary edges}, the
others are {\it interior edges}.  A vertex of valance at least $4$
will be called a {\it fork vertex}.

If $\gamma$ is a path on ${\cal D}(p/q)$, $v$ an interior vertex of
$\gamma$, and $t$ the number of triangles between the two edges of
$\gamma$ incident to $v$, then the {\em corner number\/} of $v$ on
$\gamma$ is defined as $c(\gamma, v) = t$ if the triangles are above
$\gamma$, and $c(\gamma, v) = - t$ otherwise.  Thus any minimal path
$\gamma$ from $1/0$ to $p/q$ can be written as $\gamma(c_1, ...,
c_n)$, where $c_i = c(\gamma, v_i)$ and $v_i$ is the $i$-th vertex in
the interior of $\gamma$.  Denote by $[c_1, ..., c_k]$ the partial
fraction decomposition $1/(c_1 - 1/(c_2 - ...  - 1/c_k)...)$.  Then
the rational number at the vertex $v_i$ above is $p_i/q_i = v_1 +
[c_1, ..., c_i]$.  In particular, $p/q = v_1 + [c_1, ..., c_n]$.  Note
that $v_1$ is determined by $c_1$ and $p/q$: if $m$ is the integer
such that $m<p/q<m+1$ then $v_1 = m$ if $c_1>0$, and $v_1 = m+1$ if
$c_1<0$.

\bigskip
\leavevmode

\centerline{\epsfbox{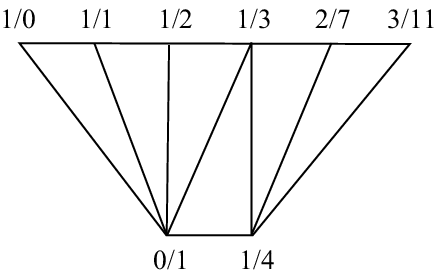}}
\bigskip
\centerline{Figure 4.3}
\bigskip

Denote by $T = T_{p/q}$ the $p/q$ rational tangle in $S[-\infty, 0]$.
It is the union of 4 vertical arcs in $S[-n,0]$ and two arcs of slope
$p/q$ on $S[-n]$ connecting the 4 endpoints of the vertical arcs on
$S[-n]$.  

Fix an orientation of $T$.  Let $\gamma$ be a path in ${\cal D}$ from
$1/0$ to $p/q$ with edges $e_1, ..., e_n$.  Let $\Sigma(e_i)$ be the
branch complex associated to the edge $e_i$ in $S[-i, -i+1]$, as
defined above.  Since the ending point of $e_i$ is the initial point
of $e_{i+1}$ and the orientations of the thick arcs are induced by
that of $T$ and hence match each other, the bottom train tracks of
$\Sigma(e_i)$ on $S[-i]$ matches the top train tracks of
$\Sigma(e_{i+1})$, so $\Sigma(\gamma) = \Sigma(e_1) \cup ...  \cup
\Sigma(e_n)$ is a branched surface with $\alpha = \Sigma(\gamma) \cap
S[0]$ a pair of train tracks of slope $1/0$ and $\beta =
\Sigma(\gamma) \cap S[-n]$ of slope $p/q$.  The bottom train tracks
are antiparallel as in Figure 4.2(c), so it can be capped off by two
copies of the trivial caps in Figure 4.1(b).  We thus obtain a
branched surface $\Sigma(\gamma)$ in the complement of $T_{p/q}$.

Let $K = K_{p/q}$ be the $p/q$ 2-bridge knot or link.  It can be
obtained from $T$ by adding two arcs of slope $1/0$ on $S[0]$
connecting the endpoints of $T$.  When the orientation of $T$ is
induced by that of $K$, the top train tracks are antiparallel, hence
they can also be capped off by trivial caps to obtain a branched
surface $\hat \Sigma(\gamma)$ in the exterior of $K_{p/q}$.
$\Sigma(\gamma)$ and $\hat \Sigma(\gamma)$ are very similar to the
branched surfaces constructed in [HT] and will be called the {\it
  Hatcher-Thurston branched surfaces\/} for $T_{p/q}$ and $K_{p/q}$,
respectively.  By Lemma 3.7 these branched surfaces are pre-laminar.
It can be shown that $\hat \Sigma(\gamma)$ is also laminar if $\gamma$
is minimal, see the proof of Theorem 5.3.

\section{Delman channels and channel surfaces}

A laminar branched surface $\Sigma$ in the exterior of a knot $K$ is
{\em persistently laminar\/} if it remains laminar in $K(r)$ for all
non-meridional slopes $r$.  To create a persistently laminar surface
we need to modify the construction of the Hatcher-Thurston branched
surface to create some meridional cusps.  The following is a
construction of Delman's channel branched surface [De].

\bigskip
\leavevmode

\centerline{\epsfbox{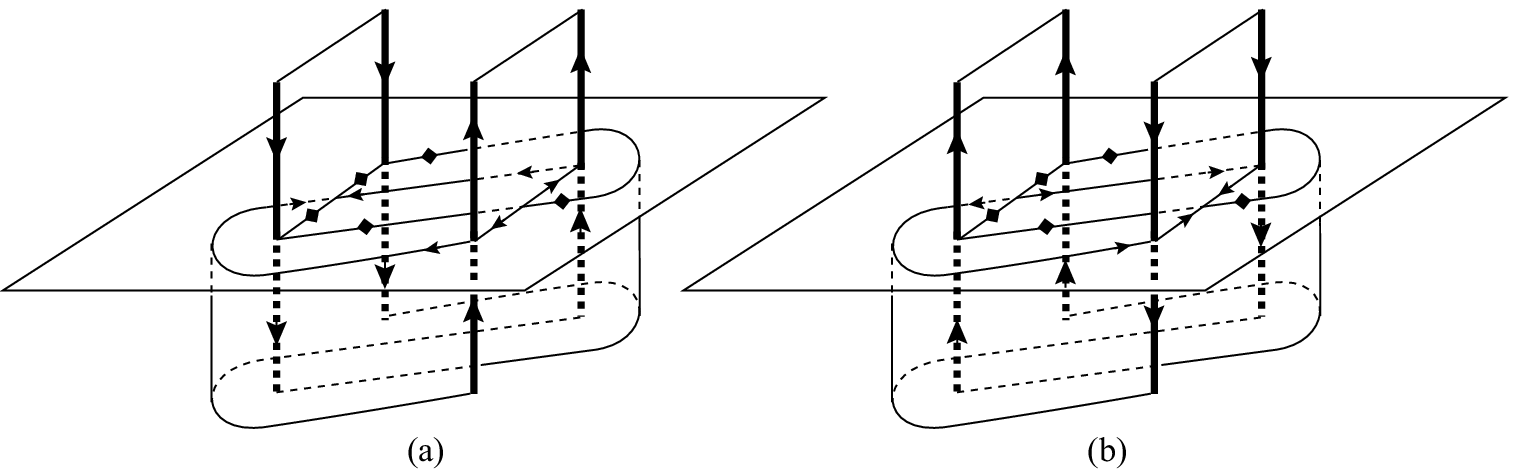}}
\bigskip
\centerline{Figure 5.1}
\bigskip

\centerline{\epsfbox{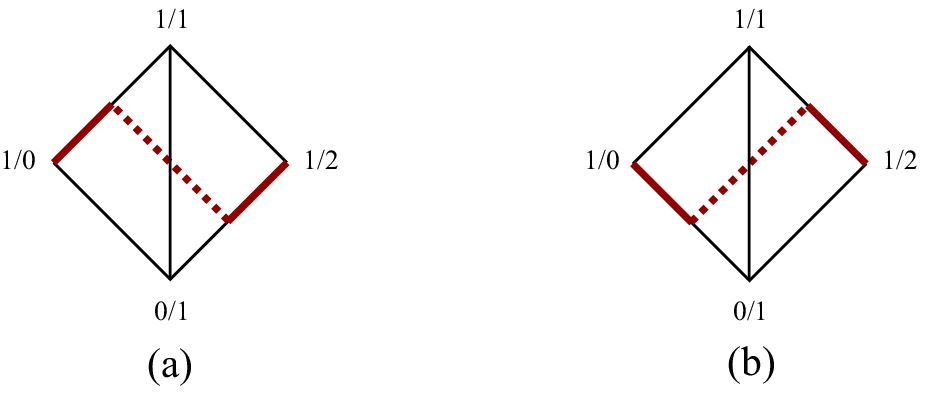}}
\bigskip
\centerline{Figure 5.2}
\bigskip

Let $\Sigma = Q \cup P \cup D$ be the 2-complex in $S[a,b]$ shown in
Figure 5.1(a), where $Q$ consists of 4 vertical tubes represented by
the 4 arcs $K$ in the figure, $P$ is a horizontal sphere, and $D$ is
the union of 4 disks of type $\alpha \times I$, where $\alpha$ is an
arc of slope $1/0$ for the two disks $D_1, D_2$ above $P$, and of
slope $1/2$ for the two disks $D_3, D_4$ below $P$.  The sink marks
are shown in Figure 5.1(a).  Let $\alpha_i = D_i \cap P$.  Unlike the
Hatcher-Thurston complex, only 2 of these arcs $\a_1, \a_3$ are
diamond edges, the other two are marked with orientation arrows.  Note
that the orientation of $\a_2, \a_4$ changes when passing across a
diamond edge, hence the two singular points in the interior of $P$ are
of type (IV) as in Definition 3.6.  The two components of $\bdd P$ on
the left are of type (I) and the two on the right are of type (II),
therefore by Lemma 3.7 these sink marks extend to a unique branched
surface structure on $\Sigma$.  The two components of $\bdd P$ of type
(II) produce two meridional cusps.

By definition the sink direction points to the right of an oriented
arc of $\bdd D_i$ when standing on $P$ on the side $D_i$ is attached
and facing to the direction of the orientation mark.  The cusp of that
arc is then on the left.  When passing across a diamond edge the cusp
continues onto the diamond edge on the side of the cusp.  Using this
fact one can check that the cusps on the top side of $P$ in Figure
5.1(a) are of slope $1/1$, while those on the bottom side of $P$ are
of slope $0/1$.

On the Hatcher-Thurston diagram $\cal D$, the branched surface in
Figure 5.1(a) is represented by a {\em Delman channel of type A} shown
in Figure 5.2(a), which is an arc starting from the vertex $1/0$,
going half way towards the vertex $1/1$, then jump to the middle of
the edge from $0/1$ to $1/2$ and finish at the vertex $1/2$.  It
reflects the properties of $\Sigma$ that the top train tracks have
slope $1/0$, the cusps on the top side of $P$ are of slope $1/1$, the
cusps on the bottom side of $P$ are of slope $0/1$, and the bottom
train tracks have slope $1/2$.

We may reverse the orientations of all the edges in Figure 5.1(a) to
obtain the one in Figure 5.1(b).  As above, these sink marks
determines a branched surface structure on the 2-complex.  By checking
the cusps one can show that the top cusps are of slope $0/1$ while the
bottom ones are of slopes $1/1$.  Therefore it corresponds to the path
in Figure 5.2(b), called a {\em Delman channel of type B}.

The branched surfaces in Figure 5.1 are called the {\em Delman channel
  surfaces}, denoted by $\Sigma(e)$ if $e$ is the corresponding Delman
channel. It is of type A or type B according to the type of $e$.  Note
that each boundary arc of $\Sigma(e)$ must connect a pair of parallel
endpoints of the arcs $K$ representing $Q$, so that the corresponding
boundary train tracks are one positive and one negative, as in Figure
4.2(a)--(b), respectively.  

By twisting the four vertical arcs around, we can isotope this surface
$\Sigma(e)$ in $S[a,b]$ to change the top slope to $r_1 = p_1/q_1$ and
the bottom slope to $r_2 = p_2/q_2$ if $\Delta(r_1, r_2) = 2$.  Hence
we can embed a Delman channel into a pair of adjacent triangles in the
Hatcher-Thurston diagram $\cal D$.  The only requirement is that the
curves $K$ must be oriented in such a way that the two train tracks of
slope $r_1$ on $S[a]$ connecting the four points of $K \cap S[a]$ must
be one positive and one negative, as in Figure 4.2(a) and (b),
respectively.  One can see that the two train tracks of slope $r_2$ on
$S[b]$ also have the same property.

Now consider $L = T_{p/q}$ or $K_{p/q}$.  It intersects the level
sphere $S[0]$ in four points.  Recall that two points of $L$ on a
level sphere are parallel if the orientations of $L$ at these two
points are both upward or both downward.  A rational number $p/q$ is
assigned a parity pair $o/e$, $e/o$ or $o/o$ if $(p,q)$ is (odd,
even), (even, odd) or (odd, odd), respectively.  Note that two arcs on
$S[0]$ with the same parity pair connect the same pair of points of
$L\cap S[0]$, so we can use the slope $1/0$, $0/1$ or $1/1$ that has
the same parity pair as $p/q$ to determine whether a curve of slope
$p/q$ connects a pair of parallel points.

\begin{defn} Let $L = T_{p/q}$ or $K_{p/q}$.  A path $\gamma$ in
  ${\cal D}(p/q)$ from $1/0$ to $p/q$ is an {\em allowable path\/} if
  it satisfies the following conditions.

(1) It is the union of edges of ${\cal D}(p/q)$ and Delman channels;

(2) it is minimal in the sense that the corner number $c(\gamma, v_i)$
defined in Section 4 satisfies $|c(\gamma, v_i)| \geq 2$ for any
interior vertex $v_i$ of $\gamma$;

(3) the label of the ending points of the Delman channels in $\gamma$
    all have the same parity pair, which is different from that of
    $p/q$, and is also different from that of $1/0$ if $L = K_{p/q}$.
\end{defn}

Write $\gamma = e_1 \cup ... \cup e_n$, where each $e_i$ is either an
edge or a Delman channel.  Note that the beginning and ending points
of a Delman channel $e_i$ have the same parity pair.  Let $r'_i,
r''_i$ be the slopes of the beginning and ending points of the $i$-th
Delman channel in $\gamma$.  Then condition (3) above implies that
these all have the same parity pair, it is different from that of
$p/q$, and if $L$ is a knot then it is also different from that of
$1/0$.  This implies that we can orient $L$ so that an arc of slope
$r'_i$ or $r''_i$ connects a pair of parallel points of $L$, so there
is a Delman channel surface $\Sigma(e_i)$ whose sink marks on $Q$
coincide with this orientation of $L$.  As in the construction of the
Hatcher-Thurston branched surfaces for minimal path, we can now
construct branched surfaces $\Sigma(\gamma)$ and $\hat \Sigma(\gamma)$
for $\gamma$ as the union of $\Sigma(e_i) \subset S[-i, -i+1]$ and
some copies of the trivial caps in Figure 4.1(b) at the top and
bottom, except that when $e_i$ is a Delman channel $\Sigma(e_i)$ is
the corresponding Delman channel surface above instead of the
Hatcher-Thurston branched surface in Figure 4.1.  The branched
surfaces $\Sigma(\gamma)$ in $S[-\infty, 0]$ and $\hat \Sigma(\gamma)$
in $S^3$ are called the {\it Delman branched surfaces\/} for the
tangle $T_{p/q}$ and the link $K_{p/q}$, respectively, corresponding
to the allowable path $\gamma$.

\begin{lemma} Let $\Sigma = \Sigma(\gamma)$ and $\hat \Sigma = \hat
  \Sigma(\gamma)$ be the Delman branched surfaces corresponding to an
  allowable path $\gamma = \gamma(c_1, ..., c_n)$ in ${\cal D}(p/q)$,
  where $c_i = c(\gamma, v_i)$ is the corner number of $\gamma$ at
  $v_i$.  Then each of $E(\hat \Sigma)$ and $E(\Sigma)$ has $n+1$
  components $Y_0, ..., Y_n$, one for each vertex $v_i$ of $\gamma$,
  such that $Y_n$ for both branched surfaces and $Y_0$ for $\hat
  \Sigma$ are 3-balls with a single cusp, and $Y_i$ is a solid torus
  with cusp winding number $|c_i|$ for $i\neq 0, n$.
\end{lemma}

\proof Topologically the exterior of $\hat \Sigma(\gamma)$ is obtained
by cutting $S^3$ along horizontal spheres $S[-i+1/2]$ for $i=1 , ...,
n$, then removing a regular neighborhood of the attaching
disks $D$ in each region.  Hence it has $n+1$ regions $Y_0, ..., Y_n$,
one for each vertex $v_i$.  It is easy to see that $Y_0$ and $Y_n$ are
3-balls with a single cusp.  For $i=1, ..., n-1$, $D \cap S[-i-1/2,\,
-i+1/2]$ is a pair of disks $D_{i1} \cup D_{i2}$, each $D_{ij}$ is a
$\alpha_j \times I$ for a curve $\alpha_j$ of slope $p_i/q_i$, hence
$Y_i$ is a solid torus.  The exterior of $\Sigma(\gamma)$ is the same,
except that $Y_0$ is now a solid torus in $S[-1/2, 0]$.

We need to determine the winding number of the cusps on $\bdd Y_i$.
By an isotopy we can deform $\hat \Sigma(\gamma)$ so that the curves
$\alpha_j$ above have slope $1/0$, and the cusps on the bottom side of
$S[-i+1/2]$ are of slope $0/1$.  The deformation changes labels of all
vertices of $\cal D$ but preserves $\Delta(r,s)$.  Whether $e_i$ is an
edge or Delman channel, the ending segment of $e_i$ now lies on the
edge from $1/0$ to $0/1$, so the initial segment of $e_{i+1}$ must be
on the edge from $0/1$ to $1/m$ for some integer $m$ because all
vertices connected to $0/1$ is of that form.  Moreover, the corner
number $c(\gamma, v_i) = m$.  It follows that the cusps on the top of
$S[-i-1/2]$ is of slope $1/m$.  It is now easy to see that the minimal
intersection number between a meridian disk of $Y_i$ and the cusp on
$\bdd Y_i$ is $|m| = |c(\gamma, v_i)| = |c_i|$.  \qed \medskip

The following theorem is an extension of a result of Delman
[De] for these branched surfaces, which has been used in [BW] to
determine small Seifert fibered surgeries on 2-bridge knots.

\begin{thm} Given an allowable path $\gamma$ from $1/0$ to $p/q$, the
  corresponding branched surfaces $\Sigma = \Sigma(\gamma)$ and $\hat
  \Sigma = \hat \Sigma(\gamma)$ are pre-laminar.  $\hat \Sigma$ is
  laminar in $E(K_{p/q})$, and is genuine if at least one vertex $v$
  on $\gamma$ has $|c(\gamma, v)| > 2$.  If $\gamma$ has $k$ Delman
  channels then $\Sigma$ and $\hat \Sigma$ have $2k$ meridional cusps.
\end{thm}

\proof The first statement follows from Lemma 5.2.  If $\gamma =
\gamma(c_1, ..., c_n)$ is an allowable path then $|c_i| \geq 2$ for
all $i$, so by Lemmas 5.2 $E(\hat \Sigma(\gamma))$ is an essential
cusped manifold, and by Lemma 2.2 $\hat \Sigma(\gamma)$ is a laminar
branched surface.  If some $|c_i| \geq 3$ then the corresponding
region $Y_i$ in Lemma 5.2  is not an $I$-bundle, hence $\hat
\Sigma(\gamma)$ is genuine.  Since each Delman channel creates two
meridional cusps, $\hat \Sigma$ and $\Sigma$ have $2k$ meridional
cusps if $\gamma$ contains $k$ Delman channels.  \qed \medskip

\section{Half channel surfaces}

We can reverse the orientations of some of the edges of the Delman
channel surface in Figure 5.2 to obtain new branched surfaces.  The
ones in Figure 6.1($a_1$)-($b_3$) are called {\em Delman half channel
  surfaces}, and the corresponding path shown in the figure are their
Delman half channels.  Some of these are constructed in [De] in a more
sophisticated way.  As before, denote by $Q$ the 4 tubes represented
by the 4 vertical arcs $K$, $P$ the horizontal punctured sphere, and
$D$ the 4 disks attached to $Q \cup P$.  We use $\Sigma(x_i) = Q\cup P
\cup D$ to denote the 2-complexes with sink marks in Figure 6.1($x_i$),
$x=a,b$ and $i=1,2,3$.  For $i=2,3$ we allow the orientations of both
of the two left edges in $\Sigma(x_i)$ be changed simultaneously, so
there are two choices of $\Sigma(x_i)$ in this case.

\bigskip
\leavevmode

\centerline{\epsfbox{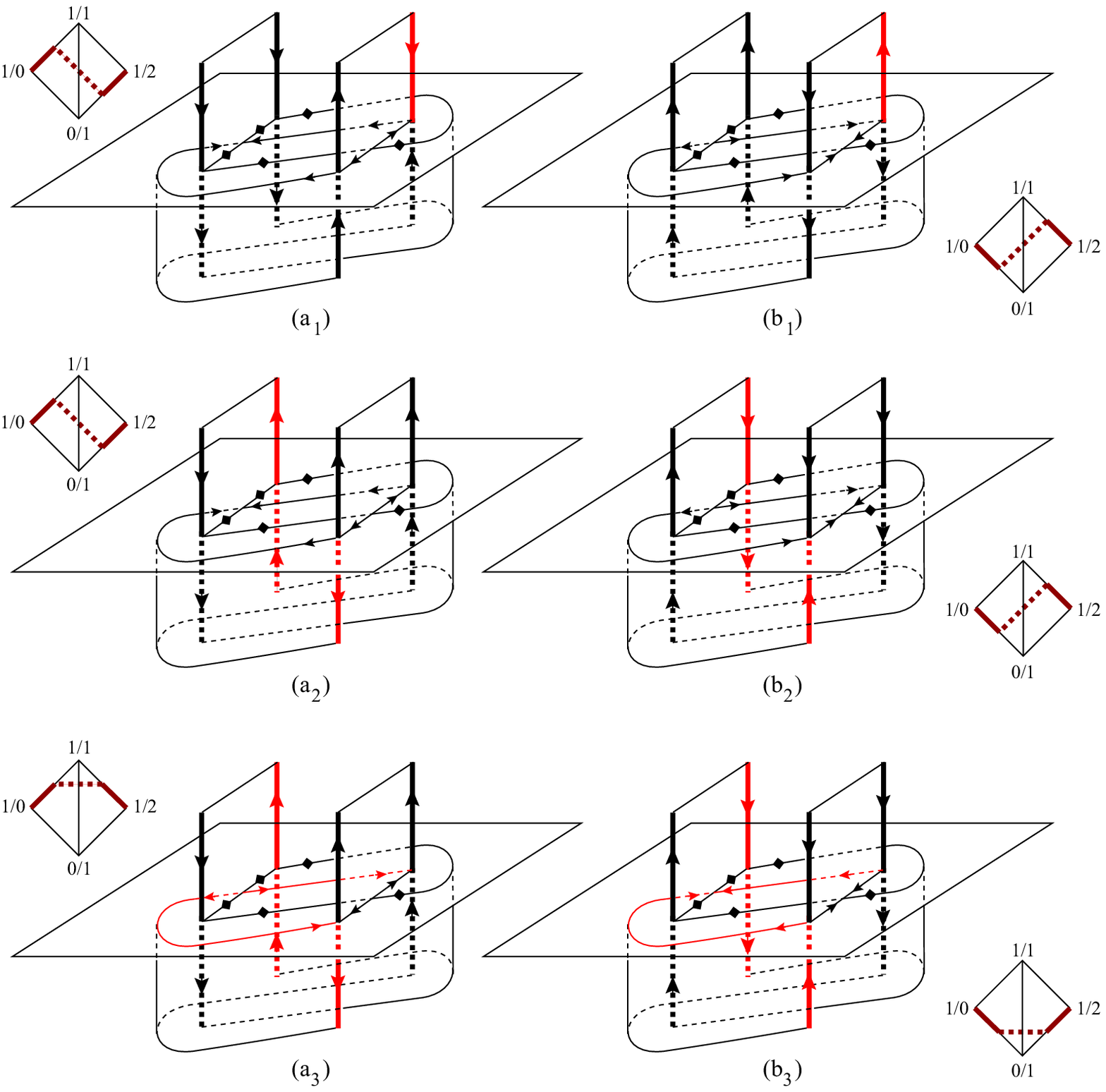}}
\bigskip
\centerline{Figure 6.1}
\bigskip

\begin{lemma} (1) $\Sigma(x_i)$ are branched surfaces.

(2) Each $\Sigma(x_i)$ has one meridional cusp.

(3) If the Delman half channel $\gamma$ representing $\Sigma(x_i)$ has
starting edge from $1/0$ to $r_1$ and ending edge from $r_2$ to $1/2$
then the cusp above $P$ has slope $r_1$ and the one below $P$ has
slope $r_2$.
\end{lemma}

\proof As for the Delman channel surfaces, (1) and (2) can be verified
using Lemma 3.7.  Note that in each case $\bdd P$ has two components
of type I, one component of type II, and one component of type III.
Therefore there is exactly one meridional cusp in each case.

(3) The cusp slopes on the two sides of $P$ are determined by the
diamond edges and the orientation sink marks on $\bdd D \cap P$.  When
$i=1,2$ the sink marks on $P$ are the same as those in Figure 5.1, so
the Delman half channel looks the same as those for the Delman
channels in Figure 5.2.  In Figure 6.1($a_3$)-($b_3$) the orientations
of $\bdd D \cap P$ have been reversed for one of the disks below $P$.
One can verify that it changes the cusp slope below $P$ to the one
indicated by the Delman half channel shown in the figure.  \qed
\medskip

Let $\gamma$ be a path in ${\cal D}(p/q)$ consisting of edges and a
single half channel $\tau$ of type $x_i$.  Denote $\Sigma(x_i)$ by
$\Sigma(\tau)$.  We can construct a 2-complex $\Sigma(\gamma)$ in the
3-ball $S[-\infty, 0]$ starting with the half channel surface
$\Sigma(\tau)$, then adding Hatcher-Thurston surfaces $\Sigma(e)$
successively for the edges $e$ before and after $\tau$.  Recall that
there are two possible choices of $\Sigma(\tau)$ if $i>1$.  We say
that $\gamma$ satisfy the {\it orientation requirement\/} if one can
choose $\Sigma(\tau)$ so that the two bottom train tracks of
$\Sigma(\tau) \cup (\cup_e \Sigma(e))$ can be capped off by trivial
caps in Figure 4.1(b) to form a branched surface, which will be
denoted by $\Sigma(\gamma)$.

We say that $\gamma$ has {\it starting slope\/} $r$ if the initial
segment of $\gamma$ is on the edge from $1/0$ to $r$.

\begin{lemma} Suppose $\gamma$ is a path in ${\cal
    D}(p/q)$ from $1/0$ to $p/q$, containing a half channel $\tau$ of
  type $x_i$.  Then $\gamma$ satisfies the orientation requirement
  unless $i=1$ and the initial point (and hence the ending point) of
  $\tau$ has the same parity pair as that of $p/q$
\end{lemma}

\proof Let $\Sigma(\tau)$ be a branched surface in $S[a,b]$
corresponding to the half channel $\tau$ of type $x_i$.  When $i>1$ we
have two choices of $\Sigma(\tau)$ and the second one is obtained from
that in Figure 6.1($x_i$) by reversing the orientations of the two
arcs on the left, hence one of them has the property that an arc of
slope $p/q$ on the bottom level sphere $S[a]$ connect a pair of
antiparallel endpoints of the 4 vertical arcs $K$ in Figure
6.1($x_i$).  The bottom train track of $\Sigma(\tau) \cup (\cup_e
\Sigma(e))$ has slope $p/q$, therefore the above implies that these
train tracks are antiparallel as in Figure 4.2(c); hence one can
cap it off using the trivial caps to obtain a branched surface
$\Sigma(\gamma)$.

Now assume $i=1$.  In this case the train tracks at the bottom of
$\Sigma(\gamma)$ are one positive and one negative, as in Figure
4.2(a)-(b).  If the endpoint of $\tau$ has different parity pair from
that of $p/q$ then an arc of slope $p/q$ connects a pair of
antiparallel edge endpoints of $K$ (with piecewise orientation induced
by sink marks), hence $\Sigma(\tau) \cup (\cup_e \Sigma(e))$ can be
capped off by trivial caps to make a branched surface.  \qed \medskip

\begin{defn} A path $\gamma$ in ${\cal D}(p/q)$ from $1/0$ to $p/q$
  consisting of edges and one half channel is a {\em semi-allowable
    path\/} if (1) it satisfies the orientation requirement above, and
  (2) it is minimal in the sense that the corner number $c( v_i,
  \gamma)$ defined in Section 4 satisfies $|c(v_i, \gamma)| \geq 2$
  for any interior vertex $v_i$ of $\gamma$.  The path $\gamma$ is
  {\em genuine\/} if $|c(v_i, \gamma)| \geq 3$ for some $i$.
\end{defn}

\begin{lemma} Let $\gamma$ be a semi-allowable path in ${\cal D}(p/q)$
  with starting slope $r$ and corner number $c_i$ at the $i$-th
  vertex.  Then the exterior of $\Sigma(\gamma)$ has $n+1$ components
  $Y_0, ..., Y_n$, one for each vertex $v_i$ of $\gamma$, such that
  $Y_n$ is a 3-ball with a single cusp, and $Y_i$ is a solid torus
  with cusp winding number $|c_i|$ for $i\neq 0, n$.  The cusps above
  the top level surface is of slope $r$.
\end{lemma}

\proof This is similar to Lemma 5.2.  We omit the details.
\qed \medskip

Note that if $\gamma$ contains a half channel and if $T_{p/q}$ is
endowed with the induced piecewise orientation, then exactly three of
the four endpoints of $T_{p/q}$ on the top level sphere $S[0]$ have
the same orientation.  We say that $\gamma$ is {\it upward\/} if the
orientations at those three points are upward, and {\it downward\/}
otherwise.  From Figure 6.1 we see that $\gamma$ is upward if and only
the half channel in it is of type $b_1$, $a_2$ or $a_3$.

\begin{prop} Suppose $0<p/q<1$.  Let ${\cal D} = {\cal D}(p/q)$.

  (1) ${\cal D}$ always has an upward semi-allowable path $\gamma$
  with starting slope $1$, and there is a genuine such $\gamma$ unless
  $p = 1$ or $q-1$.

  (2) ${\cal D}$ has a downward semi-allowable path $\gamma$ with
  starting slope $1$ unless $p = 1$, or $p=q-1$ and $q$ is even.
  ${\cal D}$ has a genuine such $\gamma$ unless $p = 1, 2$ or $q-1$.

  (3) ${\cal D}$ always has a downward semi-allowable path $\gamma$
  with starting slope $0$, and there is a genuine such $\gamma$ unless
  $p = 1$ or $q-1$.

  (4) ${\cal D}$ has an upward semi-allowable path $\gamma$ with
  starting slope $0$ unless $p = q-1$, or $p=1$ and $q$ is even.
  ${\cal D}$ has a genuine such $\gamma$ unless $p = 1, q-2$ or $q-1$.
\end{prop}

\proof (1) If $p/q \geq 1/2$, let $\gamma$ be the half channel of type
$a_2$ in Figure 6.1($a_2$) followed by a path on the lower boundary of
${\cal D}$; if $p/q < 1/2$, let $\gamma$ be the half channel of
type $a_3$ followed by a path on the upper boundary of ${\cal D}$.  By
Lemma 6.2 $\gamma$ satisfies the orientation requirement and hence is
semi-allowable.  If $p \neq 1, q-1$ then there are fork vertices on
both the upper and lower boundary of $\cal D$, so $|c(v)|>2$ for some
$v\in \gamma$, hence $\gamma$ is genuine.

(2) First assume $p/q>1/2$.  Then the vertex $1/2$ is on the lower
boundary of $\cal D$.  If $q$ is odd then we may choose $\gamma$ to be
the half channel in Figure 6.1($a_1$) followed by edge path on the
lower boundary of $\cal D$, as shown in Figure 6.2(a).  By Lemma 6.2
this is semi-allowable since $1/2$ has different parity pair from that
of $p/q$, and it is genuine if $p \neq q-1$ because in this case there
is at least one fork vertex at the lower boundary.

If $q$ is even, then $p \neq q-1$ by assumption, so the last interior
edge at $v' = 1/1$ connects it to a fork vertex $v''$ on the lower
boundary.  Let $\gamma$ be the path in Figure 6.2(b)-(c) according to
whether $v''$ has only two interior edges or more, where the half
channel is of type $a_1$.  Since $1/1$ has parity pair different from
that of $p/q$, by Lemma 6.2 $\gamma$ is semi-allowable.  It is genuine
unless $|c(v', \gamma)| = 2$ (so $\gamma$ is as shown in Figure
6.2(c)), and there is no fork vertex on the upper boundary other than
$1/1$.  However, in this case $v'' = 1/2$, and it has an edge
connected to $p/q$, so $q$ cannot be even, a contradiction.  Therefore
$\gamma$ is genuine.

Now assume $p/q<1/2$, so the vertex $1/2$ is on the upper boundary of
$\cal D$.  Let $v'$ be the first fork vertex on the upper boundary,
which exists since $p \neq 1$.  Let $\gamma_1$ and $\gamma_2$ be the
paths in Figure 6.2(d)-(e), respectively.  Each $\gamma_i$ is the
union of a single half channel $\tau$ of type $a_1$ and some boundary
edges.  It is semi-allowable unless the ending point of $\tau_i$ has
the same parity pair as that of $p/q$.  Since the ending points of
$\tau_1$ and $\tau_2$ are connected by an edge of $\cal D$, they have
different parity pair.  It follows that at least one of the $\gamma_i$
is semi-allowable.  It remains to show that $\gamma_i$ can be chosen
to be genuine if $p\neq 1, 2, q-1$.

If $\gamma_1$ is semi-allowable but not genuine, then $v'$ is the only
fork vertex on the upper boundary, and it has only two interior edges,
so $v'$, $p/q$ and all the vertices between them on the upper boundary
connect to the same vertex on the lower boundary, which must be the
vertex $v''$ in Figure 6.2(d) because there is no edge connecting $v'$
to any vertex between $v''$ and $p/q$ on the lower boundary.
Therefore $v''$ and $p/q$ are connected, so they have different parity
pairs, which implies that $\gamma_2$ is also semi-allowable.  Hence
either $|c(v'', \gamma_2)| \geq 3$ and we are done, or there is no
vertices between $v'$ and $p/q$, in which case $p/q = 1/(a+1/2) =
2/(2a+1)$, hence $p = 2$ and the result follows.

If $\gamma_1$ is not semi-allowable then $\gamma_2$ is semi-allowable,
so if $\gamma_2$ is not genuine then all vertices between $v''$ and
$p/q$ (including the two endpoints) would have edge connected to $v'$,
which implies $v'$ is connected to $p/q$ and hence has different
parity pair from that of $p/q$, so $\gamma_1$ is also semi-allowable,
which is a contradiction.

The proofs of (3) and (4) are similar, using the half channels of types
$b_i$ instead.  \qed \medskip

\bigskip
\leavevmode

\centerline{\epsfbox{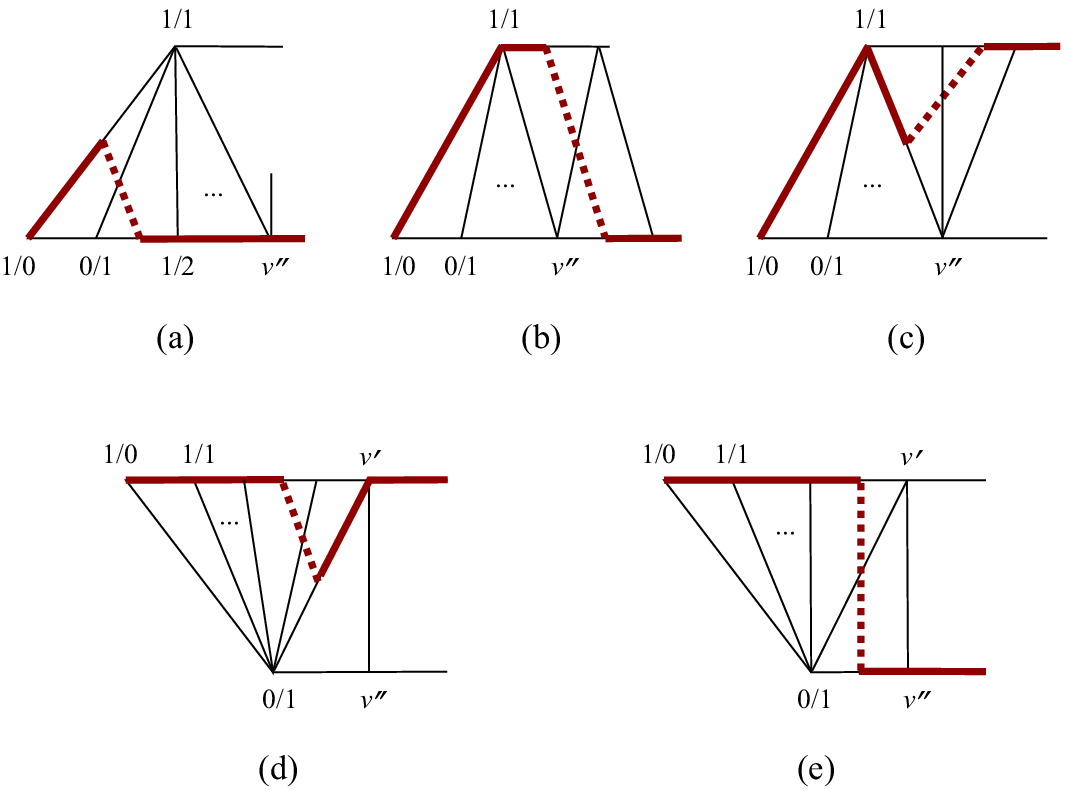}}
\bigskip
\centerline{Figure 6.2}
\bigskip

\begin{thm} Let $K$ be a non 2-bridge Montesinos knot.  Then $K$ has a
  persistently laminar branched surface in its complement unless it is
  equivalent to $K(1/q_1,\, 1/q_2,\, 1/q_3,\, -1)$, where $q_i$ are
  positive integers, and $q_1$ is even.
\end{thm}

\proof If the length of $K$ is at least 4 then by [Wu1] its exterior
of $K$ contains a closed essential surface which remains
incompressible after all surgery on $K$.  Hence we may assume that $K$
is of length 3.  Write $K = K(r_1, r_2, r_3, n)$, where $0<r_i =
p_i/q_i<1$ and $n$ is an integer.  By taking the mirror image if
necessary we may assume that $n \geq -1$.  We may assume that $q_2,
q_3$ are odd.

First assume $n\geq 0$.  By Proposition 6.5(1) there is an upward
semi-allowable path $\gamma_1$ in ${\cal D}(r_1)$ with starting slope
$1$, and there is a downward semi-allowable path $\gamma_2$ in ${\cal
  D}(r_2)$ with starting slope $0$.  Let $\gamma_3$ be the upper
boundary path in ${\cal D}(r_3)$, which has starting slope $1$.  Let
$\Sigma(\gamma_i)$ be the branched constructed above for $i=1,2$, and
let $\Sigma(\gamma_3)$ be the Hatcher-Thurston branched surface, with
the arcs of $T(r_3)$ oriented so that the two left edges points
outward.  (This is possible since $q_3$ is odd.)  Inserting these into
a band with $n$ half twist produces a branched surface $\Sigma$ in the
knot exterior as shown in Figure 6.3.  By Lemma 3.7 $\Sigma$ is
pre-laminar.  Since $\gamma_1$ and $\gamma_3$ has starting slope $1$
and $\gamma_2$ has starting slope $0$, the component of $E(\Sigma)$
outside of the tangles is the complement of a band with $n+2$ half
twists, and it has the boundary of the band as its cusps, so it is an
essential cusped manifold.  Any other component of $E(\Sigma)$ is a
component of the exterior of some $\Sigma(\gamma_i)$ in the
corresponding tangle space, which by Lemma 6.4 is a solid torus with
cusp winding number at least 2.  Hence $\Sigma$ is laminar by Lemma
2.2.  Since it has two meridional cusps, it remains laminar after all
Dehn surgery.

\bigskip
\leavevmode

\centerline{\epsfbox{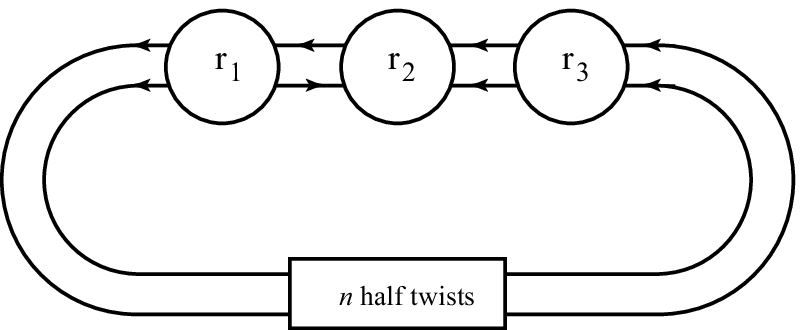}}
\bigskip
\centerline{Figure 6.3}
\bigskip

Now assume $n = -1$.  If for some $i$, $q_i$ is odd and $p_i \neq 1$,
or $q_i$ is even and $p_i \neq 1, q_i -1$, then by Proposition  6.5(2) we
have a downward semi-allowable path $\gamma_i$ with starting slope
$1$.  We can then choose $\gamma_j$ ($j=1$ if $i\neq 1$) to be an
upward semi-allowable path with starting slope $1$ and $\gamma_k$ ($k
\neq i,j$) the upper boundary path, and construct $\Sigma$ as above.
Since each of these paths has starting slope $1$, the outside
component of $\Sigma$ has $n+3 = 2$ half twists and hence for the same
reason as above, $\Sigma$ is a persistently laminar branched surface.

If $q_1$ is also odd then the above shows that $p_i = 1$ for all $i$,
in which case $K(1/p_1,\, 1/p_2,\, 1/p_3,\, -1)$ is a link of two
components, contradicting the assumption.  Therefore we may assume
$q_1$ is even.  By the above it remains to consider the case that $p_1
= q_1 - 1$, and $p_2 = p_3 = 1$.  If $q_1=2$ then $p_1 = 1$ and we are
done, so we may assume $q_1 \geq 4$.  Here the construction is
different.  In this case $\gamma_1$ is a path with two half channels,
as shown in Figure 6.4(a).  The half channel at the end is of type
$b_2$.  The one at the beginning is a new half channel with its
corresponding branched surface shown in Figure 6.4(b), which is
obtained from the one in Figure 6.1($a_2$) by rotating upside down
along the horizontal axis from left to right, and then twisting the
strings so that the top train track has slope $1/0$.  One can also
check directly from Figure 6.4(b) that it is a branched surface with a
single meridional cusp, the cusp on the top side of $P$ has slope $1$
and the one on the bottom side of $P$ has slope $0$.  Now the
orientations of the arcs at the bottom level of this surface match the
orientation of the top arcs of one of the surfaces of type $b_2$ (the
one with orientations of both left arcs of that in Figure 6.1($b_2$)
reversed).  We can then add some Hatcher-Thurston surfaces
corresponding to the edges between these two half channels if
necessary to obtain a branched surface $\Sigma(\gamma_1)$
corresponding to the path $\gamma_1$ of ${\cal D}(r_1)$ in Figure
6.4(a).  The other two tangles are of type $T_{1/q_i}$ with $q_i$ odd,
so we can orient them as shown in Figure 6.4(c).  Let $\gamma_i$ be
the upper boundary path in ${\cal D}(1/q_i)$, which has starting slope
$1$, and let $\Sigma(\gamma_i)$ be the corresponding Hatcher-Thurston
branched surface.  Inserting these into the band with $n=-1$ half
twist as in Figure 6.4(c), we obtain the branched surface $\Sigma$.
For the same reason as before, $\Sigma$ is a persistently laminar
branched surface.

We note that in the construction above we may reverse the orientations
of both left edges in the half channel surfaces in Figure 6.4(b) and
Figure 6.1($b_2$), so the outside orientations of $\Sigma(\gamma_1)$
is obtained from that in Figure 6.4(c) by reversing the orientations
of the two arcs on the left of the tangle, as shown in Figure 6.4(d).
The orientations would not match that of the right endpoints of the
third tangle, but they do when there is no half twist or an even
number of half twists at the bottom, as shown in the figure for $n=0$.
Therefore this modification produces a persistently laminar branched
surface $\Sigma$ when $n \geq 0$ is even.  This modification is needed
in the proof of Theorem 6.7.
\qed \medskip

\bigskip
\leavevmode

\centerline{\epsfbox{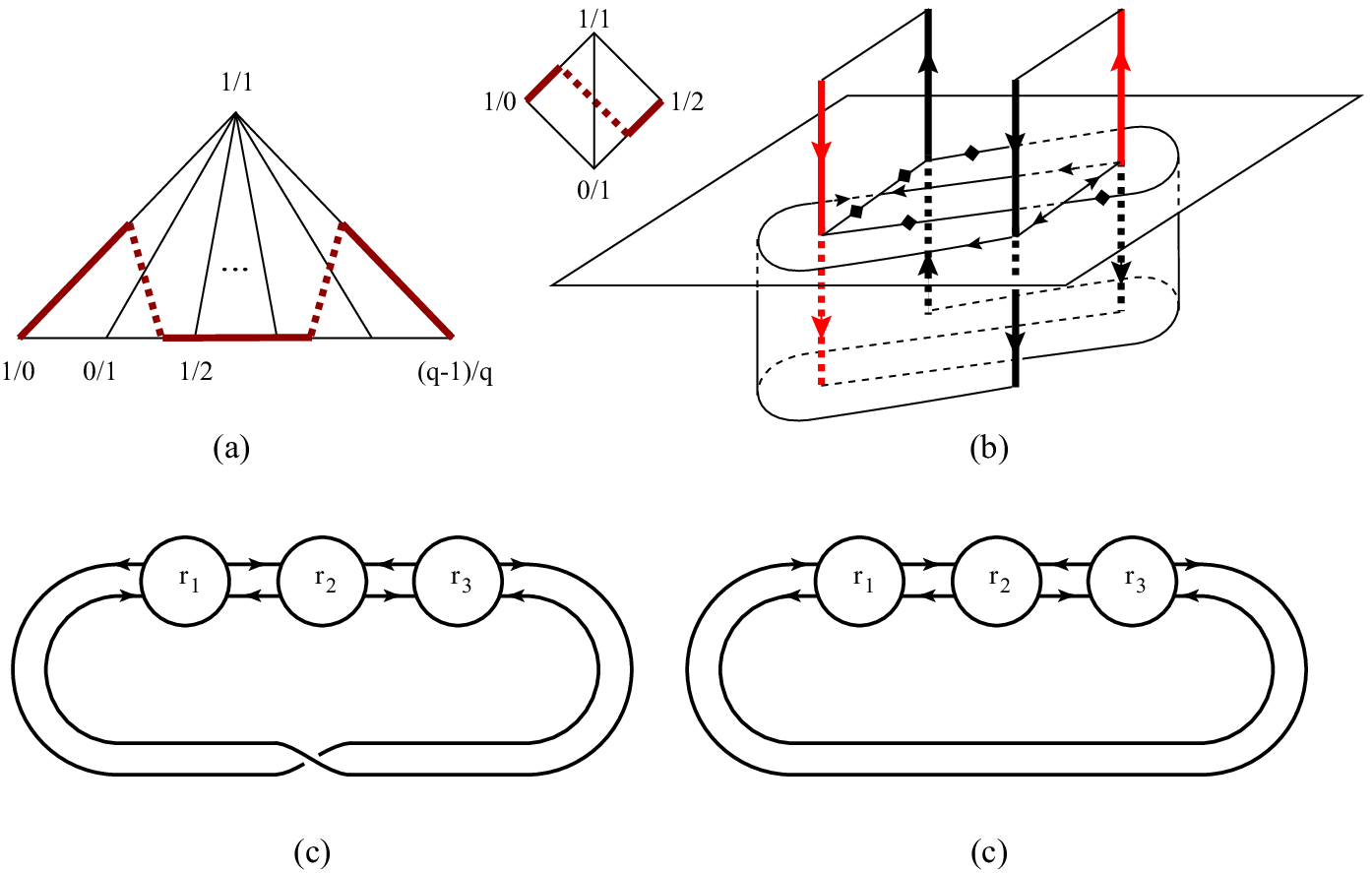}}
\bigskip
\centerline{Figure 6.4}
\bigskip

\begin{thm} Let $K$ be a Montesinos knot of length 3.  Then $K$ has a
  {\em genuine\/} persistently laminar branched surface in its
  complement unless $K$ is equivalent to $K(1/q_1,\, 1/q_2,\,
  p_3/q_3,\, n)$, such that either

  (1) $n = 0$, $q_i\geq 2$, and $p_3 = 1$; or

  (2) $n =-1$, $q_i \geq 2$, and $p_3 = 1, 2$ or $q_3 - 1$.
\end{thm}

\proof As before, we may assume that $K = K(r_1, r_2, r_3, n)$, where
$0<r_i = p_i/q_i < 1$ and $n \geq -1$.  If $n>0$, or if $n=0$ and some
$p_i \neq 1$ then we can construct $\gamma_i$ and $\Sigma$ as in the
proof of Theorem 6.6 for the case $n=-1$.  Now the outside component
of the exterior of $\Sigma$ is the exterior of a band with at least
three half twists, hence $\Sigma$ is genuine.  (Unlike Theorem 6.6, we
cannot claim that $q_1$ is even in this case, because the proof of
Theorem 6.6 used the fact that $K(1/q_1,\, 1/q_2,\, 1/q_3,\, n)$ is a
link of two components if all $q_i$ are odd and $n=-1$, which is no
longer true when $n=0$.)

We may now assume that $n = -1$.  The result is true if $p_1 = p_2 =
p_3 = 1$, or if $p_1 = p_2 = 1$ and $p_3 = q_3 - 1$.  Thus we may
assume that $p_3 \neq 1$, and either $p_1 \neq 1$ or $p_3 \neq q_3
-1$.  Note that if $p_1 \neq 1$ then up to relabeling of the tangles
we may assume that $q_3$ is odd.  Hence either $p_3 \neq q_3 -1$ or
$q_3$ is odd.  We may now apply Proposition 6.5(2) to obtain a
downward semi-allowable path $\gamma_3$ in ${\cal D}(p_3/q_3)$ with
starting slope $1$.  If $p_j \neq 1$ for some $j=1,2$, let $\gamma_j$
be the upper boundary path in ${\cal D}(p_j/q_j)$, which is genuine in
this case, and let the other $\gamma_k$ be an upward semi-allowable
path with starting slope $1$, as in Proposition 6.5(1).  If $p_1 = p_2
= 1$ but $p_3 \neq 1, 2, q_3 -1$ then by Proposition 6.5(2) we may
choose $\gamma_3$ to be genuine.  It follows that the branched surface
$\Sigma$ corresponding to these paths is genuine unless $p_1 = p_2 =
1$ and $p_3 = 1, 2$ or $q_3 -1$.  \qed

\begin{example} {\rm The knot $10_{142}$ on the knot table of [Ro] is
    the pretzel knot $K(1/3,\, 1/3,\, -1/4)$.  It is one of the 5
    knots in Gabai's Frontier Question FQ 1.2 [Ga] that were not known
    whether all surgeries are laminar.  Theorem 6.6 shows that it is
    persistently laminar.}
\end{example}

The construction of branched surfaces in rational tangle spaces can be
used to construct persistently laminar branched surfaces in the
complement of some non-Montesinos knots.  Here is an example.  Let $L$
be a non-split oriented link.  A spanning surface $F$ of $L$ is {\it
  $\pi_1$-injective\/} if it is $\pi_1$-injective in the complement of
$L$, in which case the manifold $M$ obtained by cutting $E(L)$ along
$F$ is an essential cusped manifold with $\bdd F$ as its cusp.  Let
$\alpha$ be a proper arc in $F$.  Embed a regular neighborhood $D$ of
$\alpha$ in ${\Bbb R}^2$ so that the two arcs $L \cap \bdd D = a_1
\cup a_2$ are horizontal.  Then $\alpha$ is said to {\em connect
  parallel arcs\/} if the orientations of $a_1, a_2$ points to the
same direction; otherwise $\alpha$ {\em connects antiparallel arcs.}
For example, if $F$ is a minimal Seifert surface then it is always
$\pi_1$-injective, and $\alpha$ always connect antiparallel arcs.  Set
up a coordinate on the boundary of $B = N(\alpha)$ so that $a_1 \cup
a_2$ is a $0$-tangle and $F \cap \bdd B$ is isotopic to a $\frac
10$-tangle.  Denote by $L(F, \alpha, r)$ the knot or link obtained
from $L$ by replacing $a_1 \cup a_2$ with a rational tangle $T_{r}$.

\begin{cor} Suppose $L$ is a non-split oriented link, $F$ is a
  $\pi_1$-injective spanning surface of $L$, and $\alpha$ an arc on
  $F$.  Let $K = L(F, \alpha, 1/n)$, where $|n|>2$ is odd if $\alpha$
  connects parallel arcs, and even otherwise.  If $K$ is a knot then
  it has a persistently laminar branched surface.
\end{cor}

\proof Let $B = N(\alpha)$ and $T_{1/n} = K \cap B$, as defined above.
If $n$ is odd, let $\gamma$ be an allowable path starting with the
Delman channel in Figure 5.2(b) followed by the upper boundary path,
and if $n$ is even let $\gamma$ be the path in Figure 6.4(a), except
that the labels on the top vertex is $0/1$ and the labels at the
bottom are $1/0,\, 1/1,\, ...,\, 1/n$ (so the diagram is now ${\cal
  D}(1/n)$ reflected along a horizontal line).  In the first case let
$\Sigma(\gamma)$ be the branched surface constructed in Theorem 5.3.
In the second case let $\Sigma(\gamma)$ be the one constructed in the
proof of Theorem 6.6 and, as in the proof of Theorem 6.7, we may
reverse the orientation of the two left arcs of the tangle so that the
orientation looks like that of $T_{r_1}$ in Figure 6.4(d).  The
assumption that $n$ is odd if and only if $\alpha$ connects parallel
arcs implies that the tangle can be rotated if necessary so that the
orientation of the arcs on the boundary of the tangle defined by the
Delman channel or half channel surface match those of $L-N(\alpha)$,
so we can extend $\Sigma(\gamma)$ to a branched surface $\Sigma$ in
the complement of $K$ by adding the tubes $Q$ around $L - \Int B$ and
the surface $F - \Int B$, which is attached to $Q$ using the
orientation of $L$ as its sink marks.  By Lemma 3.7 $\Sigma$ is
pre-laminar.  As before, the components of $E(\Sigma)$ inside of the
ball $B = N(\alpha)$ are essential cusped manifolds, and the
assumption that $F$ is $\pi_1$-injective implies that the outside
component $E(\Sigma)$, which is the same as $E(L)$ cut along $F$, is
also an essential cusped manifold.  Therefore by Lemma 2.2 $\Sigma$ is
laminar.  By construction $\Sigma$ has two meridional cusps, hence it
remains laminar after all nontrivial Dehn surgery on $K$.  \qed

\section{Seifert fibered surgery on Montesinos knots}

Exceptional Dehn surgeries on arborescent knots have been determined
except for atoroidal Seifert fibered surgeries on Montesinos knots of
length 3.  The following is a result in this direction.

\begin{thm} {\rm [Wu3, Theorem 1.1]\/} Suppose $K = K(\frac{p_1}{q_1},
  \frac{p_2}{q_2}, \frac{p_3}{q_3})$ is a Montesinos knot of length
  3 and $q_i \geq 2$.  If $\frac 1{q_1-1} + \frac 1{q_2-1} + \frac
  1{q_3-1} \leq 1$ then $K$ admits no atoroidal Seifert fibered
  surgery.
\end{thm}

Theorem 6.7 can be used to strengthen this result by adding
restrictions to $p_i$.  We separate two cases.  Recall that a
Montesinos knot $K$ of length 3 is a pretzel knot of length 3 if it
can be written as $K(1/q_1,\, 1/q_2,\, 1/q_3,\, n)$ for some integers
$n$ and $|q_i|\geq 2$, and it is a genuine pretzel knot if in addition
$n$ = 0.

\begin{thm} Let $K$ be a pretzel knot of length 3.  If $K$ admits an
  atoroidal Seifert fibered surgery, then $K$ is equivalent to
  $K(\frac 1{q_1}, \frac 1{q_2}, \frac 1{q_3}, n)$ such that either
  $n=0$ and hence $K$ is a genuine pretzel knot, or $n = -1$ and
  $q_i>0$.  In either case $q_i$ satisfy $\frac 1{|q_1|-1} + \frac
  1{|q_2|-1} + \frac 1{|q_3|-1} > 1$.
\end{thm}

\proof The second part follows from Theorem 7.1, so we only need to
prove the first part.  By [Br] $K(r)$ is not an atoroidal Seifert
fibered manifold if it contains a genuine laminar branched surface,
hence by Theorem 6.7 the result holds except that we may have $K =
K(1/q_1,\, 1/q_2,\, p_3/q_3,\, -1)$, where $q_i \geq 2$ and $p_3 =
1,2$ or $q_3 - 1$.  We are done when $p_3 = 1$.  Since $K$ is assumed
to be a pretzel knot, the case $p_3 = 2$ does not occur unless $q_3 =
3$, in which case we also have $p_3 = q_3 -1$.  If $p_3 = q_3 - 1$
then $K = K(1/{q_1},\, 1/{q_2},\, (q_3-1)/q_3,\, -1) = K(1/{q_1},\,
1/{q_2}, 1/(-q_3) )$, which is a genuine pretzel
knot.  \qed \medskip

We now consider the case that $K = K(r_1, r_2, r_3)$ is not a pretzel
knot.  As in [Wu3], we use $\bar p = \bar p(p,q)$ to denote the mod $q$
inverse of $-p$ with minimal absolute value, i.e., $\bar p$ satisfies
$p\bar p \equiv -1$ mod $q$, and $2|\bar p| \leq q$.  We can combine
[Wu3, Theorem 8.2] with Theorem 6.7 to obtain the following result for
atoroidal Seifert fibered surgery on non-pretzel Montesinos knot.

\begin{thm} Let $K$ be a Montesinos knot of length 3.  If $K$ is not a
  pretzel knot and $K$ admits an atoroidal Seifert fibered surgery
  $K(r)$, then $K$ is equivalent to one of the following.

  (a) $K(-2/3,\, 1/3,\, 2/5)$;

  (b) $K(-1/2,\, 1/3,\, 2/(2a + 1)\, )$ and $a \in \{3,4,5,6\}$.

  (c) $K(-1/2,\, 1/q,\, 2/5)$ for some $q\geq 3$ odd;
\end{thm}

\proof By [Wu3, Theorem 8.2], $K$ is equivalent to one of the
following.

(1) $K(1/3,\, \pm 1/4,\, p_3/5)$, $p_3 \equiv \pm 1$ mod $5$;

(2) $K(1/3,\, \pm 1/3,\, p_3/q_3)$, $|\bar p_3| \leq 2$;

(3) $K(1/2,\, 2/5,\, p_3/q_3)$, $q_3 = 5$ or $7$;

(4) $K(1/2,\, 1/q_2,\, p_3/q_3)$, $q_2 \geq 5$ and $|\bar p_3| \leq
2$;

(5) $K(1/2,\, 1/3,\, p_3/q_3)$, $|\bar p_3| \leq 6$.

Since $K$ is not a pretzel knot, by Theorem 6.7 and [Br] we have $K =
K(1/q_1,\, 1/q_2,\, 2/q_3, -1)$ with $q_i$ positive and $q_3\geq 5$,
so it cannot be of type (1) above, and if it is of type (3) then $p_3
= 1 - q_3$, which is in conclusion (c).  We may write $q_3 = 2a + 1$,
so $K = K(1/q_1,\, 1/q_2,\, 2/(2a+1), -1)$.  Note that $a>1$ as
otherwise $K$ would be a pretzel knot.  Since $2a \equiv -1$ mod
$q_3$, we have $\bar p_3 = a$.  Therefore if $K$ is in (2) or (4) then
we have $|\bar p_3| = a = 2$, so $q_3 = 5$ and $K$ is in conclusion
(a) or (c).  Finally if $K$ is of type (5) then $2\leq a \leq 6$.
When $a=2$ $K$ is in conclusion (c), and when $a=3,4,5,6$ it is in
conclusion (b).  \qed \medskip

\section {Persistently laminar tangles}

Given a 2-string tangle $(B, T)$, we can add another 2-string tangle
$(B_1, T_1)$ to it to make it a knot or link $K = T\cup T_1$ in $S^3 =
B \cup B_1$, called an {\em extension\/} of $T$.  The gluing map
$\varphi: \bdd B \to \bdd B_1$ is an orientation reversing map,
sending a curve of slope $r$ on $\bdd B$ to a curve of slope $-r$ on
$\bdd B_1$.  The extension is $s$-nontrivial, or simply {\it
  nontrivial\/} when $s=0$, if there is no disk $D$ in $B_1$
separating the two strings of $T_1$, with $\bdd D$ a slope $s$ curve
on $\bdd B$.  It is a {\it pretzel extension\/} if $(B_1, T_1, s)$ can
be isotoped so that $s$ is a horizontal loop on $\bdd B_1$ and $T_1$
is a pair of vertical arcs.

\begin{defn} A closed branched surface $\Sigma$ in $E(T) = B - \Int
  N(T)$ is {\em persistently laminar\/} with degeneracy slope $s$,
  if it is laminar in $K(r)$ for all $s$-nontrivial extensions $K$
  of $T$ and all nontrivial slopes $r$ of $K$.  In this case
  $(B,T)$ is called a persistently laminar tangle. 
\end{defn}

Brittenham [Br] showed that the tangle $T(1/3,\, -1/3)$ is
persistently laminar.  It was proved by Youn [Yo] that the tangle
$T(1/3,\, -1/5)$ is also persistently laminar.  Using the techniques
developed above, we can now construct many more persistently laminar
branched surfaces.  See Theorem 8.5 below.

Suppose $\Sigma$ is a branched surface in $E(T)$ for some tangle
$(B,T)$.  Denote by $E_T(\Sigma) = E(T) - \Int N(\Sigma)$ and call it
the exterior of $\Sigma$ in the tangle space.  A component of
$E_T(\Sigma)$ is an {\em outside component} if it intersects $\bdd
E(T)$, otherwise it is an {\em inside components}.  If the outside
component $Y$ is a collar $\bdd E(T) \times I$ of $\bdd E(T)$ then we
say a curve on $\bdd Y - \bdd E(T)$ has slope $r$ if it is isotopic in
$\bdd E(T) \times I$ to a curve of slope $r$ on $\bdd B$.

\begin{lemma} Let $(B,T)$ be a tangle, and let $\Sigma$ be a closed
  pre-laminar branched surface in $\Int E(T)$.  If the inside
  components of $E_T(\Sigma)$ are essential cusped manifolds, and the
  outside component $Y$ of $E_T(\Sigma)$ is a collar of $\bdd E(T)$
  with a single cusp of slope $s$, then $\Sigma$ is a persistently
  laminar branched surface with degeneracy slope $s$
\end{lemma}

\proof  Let $(S^3, K) = (B,T) \cup (B_1, T_1)$ be an $s$-nontrivial
extension of $(B,T)$ and consider $\Sigma$ as a branched surface in
$S^3$.  Then the exterior of $\Sigma$ is a union of a component $X$
which contains the knot $K$, and the inside components of $\Sigma$ in
$B$, which by assumption are essential cusped manifolds.  Since the
outside component $Y$ of $E_T(\Sigma)$ is a collar of $\bdd E(T)$ with
a single cusp of slope $s$, we see that $X$ is the union of $B_1$ with
two 1-handles $V_1, V_2$ (i.e.\ a regular neighborhood of $T$ in $B$)
attached, and $K$ is the union of $T_1$ with the cores of $V_i$.  The
cusp on $\bdd Y$ becomes a cusp $\gamma$ on $\bdd X$ which is of slope
$-s$ on $\bdd B_1$.  By assumption $\gamma$ does not bound a disk in
$B_1 - T_1$, which implies that it does not bound a disk in $X-K$.

By assumption $\Sigma$ is pre-laminar and the inside components of
$\Sigma$ in $B$ are essential cusped manifold, therefore to prove that
$\Sigma$ is a laminar branched surface in $K(r)$, by Lemma 2.2 it
suffices to show that $X(r)$, the manifold obtained from $X$ by $r$
surgery on $K\subset X$, is an essential cusped manifold.

The frontier of $N(T_1)$ is a pair of annuli $A = A_1 \cup A_2$ which
cut $X$ into $E(T_1)$ and a solid torus $V$ with $K$ as its core.
Since $K$ is a nontrivial extension, the curve $\gamma$ does not bound
disk in $V-K$, hence the surface $\bdd E(T_1) - \gamma \cup A$, which
is a union of two pairs of pants, is incompressible in $E(T_1)$, so
$E(T_1)$ is an essential cusped manifold when considering $N(\gamma)
\cup A$ as a vertical surface.  After performing Dehn surgery on $K$
the solid torus $V$ becomes another solid torus $V(r)$ with $A$ as a
pair of non-meridional annuli, hence it is also an essential cusped
manifold when considering $A$ as vertical surface.  It is easy to show
that the union of essential cusped manifolds along vertical annuli
is still an essential cusped manifold.  Hence $X(r) = E(T_1) \cup
V(r)$ is an essential cusped manifold.  
\qed \medskip

\bigskip
\leavevmode

\centerline{\epsfbox{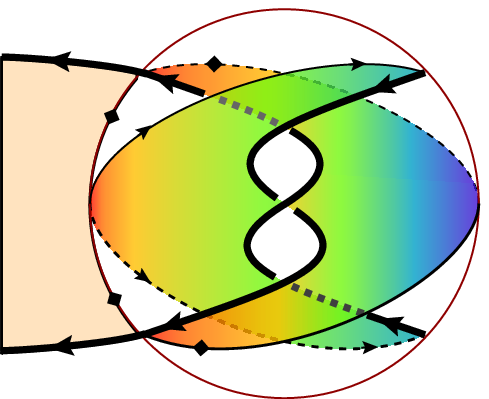}}
\bigskip
\centerline{Figure 8.1}
\bigskip

\begin{example} {\rm
Consider the 2-complex $\Sigma = Q \cup P \cup D$ in Figure 8.1,
where $Q$ is a pair of tubes represented by the thick arcs, $P$ is a
punctured sphere with two boundary components $a_1, a_2$ glued to two
boundary components of $Q$ and the other two boundary components $b_1,
b_2$ attached to the interior of $Q$, and $D$ is a  union of two
disks, one on each side of $P$.  The sink marks are shown 
in the Figure.  Note that there is only one double points of $D$,
which is of type IV, as shown in Figure 3.3(a).  The other singular
points are on the two boundary components $b_1, b_2$ of $P$, which
are of type (I) as in Definition 3.6.  Therefore by Lemma 3.7 these sink
marks can be extended to a branched surface structure for $\Sigma$. 

The train track of $\bdd \Sigma$ is positively oriented as shown in
Figure 4.2(a).  One can check that the inside component of $E(\Sigma)$
is a solid torus with cusp winding number 2.  The union of $\Sigma$
and its mirror image form a closed branched 
surface $\hat \Sigma$ in the exterior of $T(1/3,\, -1/3)$.  By Lemma
3.7 $\hat \Sigma$ is pre-laminar, and it is easy to see that it
satisfies the conditions of Lemma 8.2.  Therefore by that lemma $\hat
\Sigma$ is a persistently laminar branched surface in $T(1/3,\,
-1/3)$.  It is slightly different from the one in [Br].

It is important to note that the orientations of the sink marks make big
difference here. 
One can reverse the orientations of all oriented sink marks to obtain
another branched surface with positively oriented boundary train
track, but the inside component has cusp winding number 1 and hence is
{\em not\/} an essential cusped manifold because it has a monogon,
which is why $T(1/3,\, 1/3)$ is not a persistently laminar tangle.

The above construction can be easily generalized.  One can add more
vertical crossings to obtain a branched surface $\Sigma(1/q)$ for
$T(1/q)$ when $q>3$ is odd (but not if $q$ is even because the
orientations will not match).  One can show that the inside component
of the exterior of $\Sigma(1/q)$ is an essential cusped manifold,
hence the union of $\Sigma(1/q_1)$ and the mirror image of
$\Sigma(1/q_2)$ is a persistently laminar branched surface for
$T(1/q_1,\, -1/q_2)$ if $q_1, q_2$ are odd and $q_i \geq 3$.  We can
also replace the crossing in the middle by a horizontal band with
several crossings or other surfaces to create persistently laminar
branched surfaces for more complicated tangles.  }
\end{example}

To get more general results, we need to modify Delman's channel
surfaces.  Let $e$ be the Delman channel of type A.  Let $\Sigma'(e)$
be the the branched surface obtained from the one in Figure 5.1(a) by
deleting the upper half of the two tubes connected to the meridional
cusps, as well as the disk attached to it.  See Figure 8.2(a).  For
the same reason as before, this is a pre-laminar branched surface, and
the part below the horizontal sphere is the same as that of
$\Sigma(e)$.  Similarly, if $e$ is a Delman channel of type B then
denote by $\Sigma'(e)$ the branched surface shown in Figure 8.2(b).

Given an allowable path $\gamma$ in ${\cal D}(p/q)$ {\em starting with
  a Delman channel $e$}, the branched surface $\Sigma(\gamma)$ starts
with $\Sigma(e)$, which can be replaced by $\Sigma'(e)$ above to
obtain a new branched surface $\Sigma'(\gamma)$.  All the inside
components of $E(\Sigma'(\gamma))$ are the same as those of
$E(\Sigma(\gamma))$.  The boundary train track of $\Sigma'(\gamma)$
has only one component, which is negatively oriented if $e$ is of type
A, and positively oriented if $e$ is of type B.  The cusp on the
outside of the top horizontal surface has slope $1/2$ in both cases.

\bigskip
\leavevmode

\centerline{\epsfbox{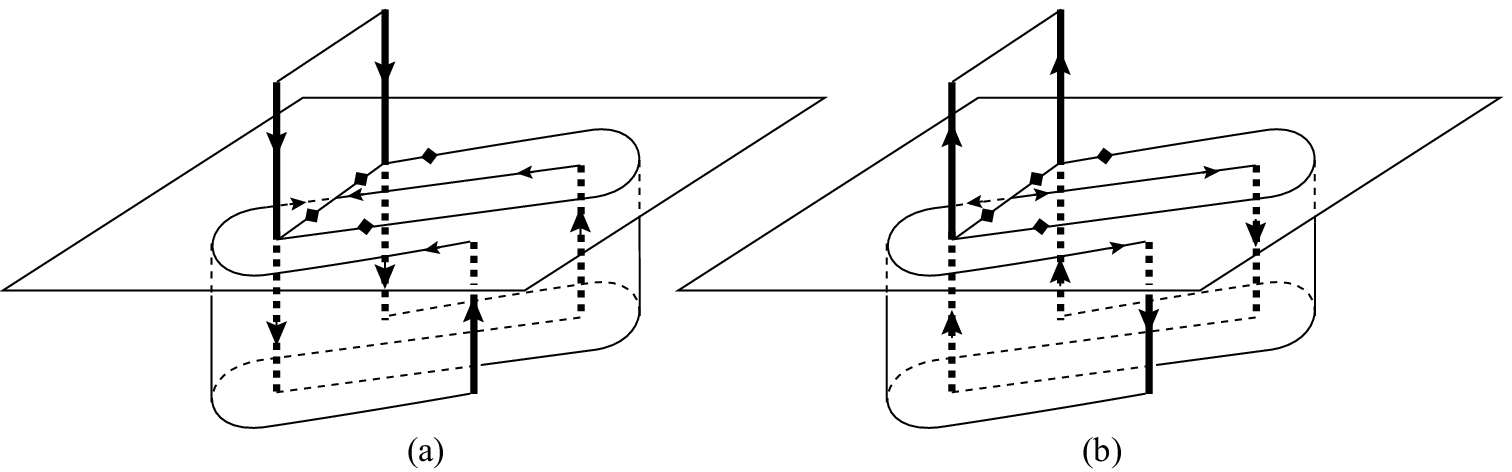}}
\bigskip
\centerline{Figure 8.2}
\bigskip

\begin{lemma} Suppose $r = p/q$, $q\geq 3$ is odd and $0<r<1$.

  (1) If $\frac 13 <r<1$ then ${\cal D}(r)$ has an allowable path
  $\gamma$ starting with a type A Delman channel.  

  (2) If $0<r<\frac 23$ then ${\cal D}(r)$ has an allowable path
  $\gamma$ starting with a type B Delman channel. 
\end{lemma}

\proof (1) If $\frac 12 < r < 1$ then ${\cal D}(r)$ is as shown in
Figure 8.3(a), in which case we can choose $\gamma$ to start with the
type A channel followed by the lower boundary path from $\frac 12$ to
$\frac pq$, as shown in the figure.  Similarly if $\frac 13<r \frac
12$ then ${\cal D}(r)$ and $\gamma$ are as shown in Figure 8.3(b).

(2) This is similar.  In this case the diagram ${\cal D}(r)$ is
obtained from that in Figure 8.3 by reflecting along a horizontal
line, then changing the label $p_i/q_i$ to $(q_i - p_i)/q_i$.  The
image of the path in the corresponding figure gives the $\gamma$
required.  \qed \medskip

\bigskip
\leavevmode

\centerline{\epsfbox{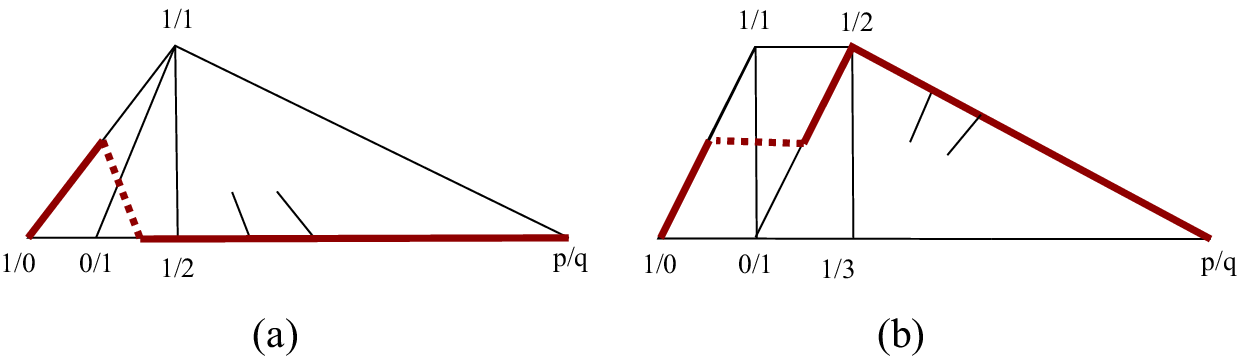}}
\bigskip
\centerline{Figure 8.3}
\bigskip

Any Montesinos tangle of length 2 can be written as $T(r_1, -r_2)$,
where $r_i = p_i/q_i$.  We can isotope it to one with $1\leq p_i
<q_i$.  (This may change the degeneracy slope $s$.)  Since $T(r_1,
-r_2)$ is homeomorphic to $T(r_2,\, -(1-r_1))$, we may also assume
that $r_1 + r_2 \leq 1$.  The following shows that many of these are
persistent tangles if both $q_i$ are odd.

\begin{thm} If $0<r_i = p_i/q_i <\frac 23$ and $q_i$ are odd then $T =
  T(r_1, \, -r_2)$ is persistently laminar.
\end{thm}

\proof Choose $\gamma_i$ to start with a type B Delman channel.  Then
$\Sigma_i = \Sigma'(\sigma_i)$ is a branched surface with positively
oriented boundary train track.  The mirror image $-\Sigma_2$ of
$\Sigma_2$ is a branched surface for $T(-r_2)$ with negatively
oriented boundary train track.  Let $\Sigma$ be obtained by gluing
$\Sigma_1$ to $-\Sigma_2$ along their boundary.  The cusps on
$\Sigma_i$ are of slope $1/2$, so $\Sigma$ has a cusp of slope $1/2 -
1/2 = 0$ on the outside.  By Lemma 8.2 $\Sigma$ is a persistent
laminar branched surface with degeneracy slope $0$.  \qed

\begin{example} {\rm Suppose $r_i = p_i/q_i$, $1\leq p_i \leq q_i$
    and $q_i$ are odd.  Then the following tangles are all
    persistently laminar.  (i) $T(1/q_1, -1/q_2)$; (ii) $T(r_1,
    -r_1)$; and (iii) $T(r_1, -r_2))$ with $r_2 \in (\frac 13, \frac
    23)$.  This follows from Theorem 8.5 for those in (i), as well as
    those in (ii) and (iii) when $r_1 < \frac 23$.  If $r_1 \geq \frac
    23$, we have $T(r_1,\, -r_1) = T(r_1-1,\, 1-r_1) = T(1-r_1,\,
    -(1-r_1))$, and $0 < 1-r_1 < \frac 23$, so (ii) also
    follows from Theorem 8.5.  Similarly if $r_1\geq \frac 23$ in case
    (iii), we have $T(r_1,\, -r_2) = T(r_1-1,\, 1-r_2) = T(1-r_2,\,
    -(1-r_1))$ and $1-r_i \leq \frac 23$, hence the result follows.}
\end{example}

\bigskip

\noindent
Department of Mathematics,  University of Iowa,  Iowa City, IA 52242
\\
Email: {\it wu@math.uiowa.edu}

\enddocument